\documentclass[11pt,a4paper]{article}
\usepackage[a4paper,margin=1in]{geometry}
\usepackage{amsmath,amssymb,amsthm,mathtools,bm}
\usepackage{graphicx}
\usepackage{subcaption}
\usepackage{booktabs}
\usepackage{enumitem}
\usepackage{microtype}
\usepackage{lmodern}
\usepackage[T1]{fontenc}
\usepackage[colorlinks=true,linkcolor=blue,citecolor=blue,urlcolor=blue]{hyperref}
\usepackage[nameinlink,noabbrev]{cleveref}
\hypersetup{
  pdftitle={Ruin Probabilities under Stochastic Premium Income and Clustered Claims},
  pdfsubject={A revised manuscript on ruin models with stochastic premiums and clustered claims},
  pdfkeywords={ruin theory, stochastic premium income, multiple iterated Poisson process, clustered claims, adjustment coefficient, Laplace transform}
}
\graphicspath{{figures/}}
\allowdisplaybreaks
\numberwithin{equation}{section}

\newtheorem{theorem}{Theorem}[section]
\newtheorem{corollary}[theorem]{Corollary}
\newtheorem{proposition}[theorem]{Proposition}

\theoremstyle{definition}

\newtheorem{example}[theorem]{Example}
\theoremstyle{remark}
\newtheorem{remark}[theorem]{Remark}

\newcommand{\E}{\mathbb E}
\newcommand{\Pbb}{\mathbb P}

\newcommand{\cL}{\mathcal L}

\newcommand{\Ftilde}{\widetilde F}
\newcommand{\Gtilde}{\widetilde G}
\newcommand{\lambdatilde}{\widetilde\lambda}
\newcommand{\mutilde}{\widetilde\mu}
\newcommand{\xitilde}{\widetilde\xi}
\newcommand{\Ntilde}{\widetilde N}
\newcommand{\Vtilde}{\widetilde V}
\newcommand{\qtilde}{\widetilde q}
\newcommand{\gtilde}{\widetilde g}

\newcommand{\hatF}{\widehat F}

\newcommand{\hatH}{\widehat H}

\title{Ruin Probabilities under Stochastic Premium Income and Clustered Claims}

\author{
Dongdong Hu \thanks{Yiwu Industrial \& Commercial College, Yiwu, China. Email: \underline{hudongdong@ywicc.edu.cn}} \and Hasanjan Sayit \thanks{Department of Financial and Actuarial Mathematics, Xi'an Jiaotong Liverpool University, Suzhou, China. Email: \underline{Hasanjan.Sayit@xjtlu.edu.cn}}}

\date{}

\begin{document}
\maketitle

\begin{abstract}
We investigate ultimate ruin and survival probabilities in insurance risk models with stochastic premium income and clustered claim arrivals. Cluster sizes are represented by multiple iterated Poisson processes (MIPPs), which provide a flexible family of compound-Poisson counting processes. We first derive a renewal-type equation for the survival probability when both premium and claim counts are MIPP-driven, and establish an exponential lower bound through an adjustment-coefficient martingale. We then study the special cases of ordinary Poisson premium arrivals and compound-Poisson premium counts, including transform formulae under exponential aggregate premium assumptions. For models that combine a deterministic premium rate with stochastic premium payments, we obtain an integro-differential equation, a corresponding adjustment-coefficient bound, and a Laplace-transform representation. Finally, we examine three classes of models with bounded claim-cluster sizes and derive explicit survival functions for selected exponential-claim specifications. The examples illustrate how cluster-size distributions and the composition of premium income affect survival probabilities when the mean surplus drift is held fixed.
\end{abstract}

\noindent\textbf{Keywords:} ruin theory; stochastic premium income; multiple iterated Poisson process; clustered claims; adjustment coefficient; Laplace transform.

\section{Introduction}

Ruin theory provides a probabilistic framework for assessing the financial strength and solvency of insurance portfolios. A central quantity is the probability that an insurer's surplus becomes negative at some time in the future. In the classical Cram\'er--Lundberg model, premiums accumulate at a deterministic constant rate, while claims arrive according to a homogeneous Poisson process and have independent, identically distributed (i.i.d.) sizes. This model has generated a large body of actuarial theory, but its assumptions may be restrictive when premium receipts are irregular or claims occur in groups. In particular, a deterministic premium stream cannot represent variation in the timing of premium payments, and a simple Poisson claim count does not directly reflect the possibility that one event produces several related claims.

One way to relax the classical assumptions is to model premium income as a stochastic process. Premium payments may arrive at random times and in random amounts, so the insurer's surplus increases through a compound counting process rather than through deterministic linear accumulation. Boikov \cite{Boikov2003} studied a Cram\'er--Lundberg model with a stochastic premium process and derived an integral equation for the survival probability together with exponential bounds based on an adjustment coefficient. Related questions involving stochastic income and discounted penalty functions have been considered by Labb\'e and Sendova \cite{LabbeSendova2009}, Zhao and Yin \cite{ZhaoYin2012}, and other authors. Miao and Sendova \cite{MiaoSendova2024} discuss the implications of premium-process uncertainty for solvency assessment, while Wang, Landriault and Li \cite{WangLandriaultLi2021} study risk processes with generalized stochastic income. See also other relevant papers \cite{DibuJacob2022, XieZhou2011, GuanHu2022}.

A second modelling issue is the occurrence of claims in clusters. Catastrophes, natural disasters, and large-scale accidents may trigger multiple claims within a short period. Compound-Poisson models with random cluster sizes provide a natural representation of such events. Ramsay \cite{Ramsay1986}, for example, studied ruin probabilities in a compound Poisson--Markov chain model in which the number of claims within a group is governed by a Markov chain. More generally, dependence and clustering in claim counts have motivated several extensions of collective risk models; see, among others, Albrecher and Boxma \cite{AlbrecherBoxma2004} and Albrecher, Constantinescu and Loisel \cite{AlbrecherConstantinescuLoisel2011}. For a comprehensive treatment of ruin theory, see Asmussen and Albrecher \cite{AsmussenAlbrecher2010}.

The multiple iterated Poisson process (MIPP) offers a convenient class of counting processes for describing clustered arrivals. It is obtained through repeated Poisson compounding and can have jumps of different sizes. Buchak and Sakhno \cite{BuchakSakhno2017,BuchakSakhno2018} studied compositions of Poisson processes and related subordinated constructions. More recently, Hu, Rachev, Sayit, Yang and Yildirim \cite{HuEtAl2026} applied iterated Poisson processes to catastrophic risk modelling and derived distributional properties relevant to ruin theory.

The present paper brings stochastic premium income and clustered claims into a common ruin-theoretic framework. We consider (i) models in which both premium and claim counts are MIPP-driven; (ii) models with a constant premium rate in addition to stochastic premium income; and (iii) several special cases in which claim-cluster sizes have finite support. For these models we derive survival-probability equations, adjustment-coefficient martingales, and Laplace-transform formulae under exponential assumptions. We also give explicit survival functions in selected examples and compare different premium and claim specifications under matched mean-drift conditions.

The paper is organized as follows. Section~\ref{sec:mipp} introduces the MIPP model and derives the general survival equation and adjustment-coefficient bound. Section~\ref{sec:constant-plus-stochastic} studies models combining deterministic and stochastic premium income. Section~\ref{sec:examples} develops three finite-cluster examples and presents explicit formulae and numerical illustrations. Concluding observations are given in Section~\ref{sec:conclusion}.

\section{Stochastic premium income driven by multiple iterated Poisson processes}
\label{sec:mipp}

We begin with a model in which both premium income and claim payments may arrive in clusters. The cluster-arrival mechanisms are represented by independent MIPPs, while the individual premium and claim amounts are assumed to be i.i.d. within each sequence.

\subsection*{MIPP construction and notation}
Let $N$ be a Poisson process with rate $\lambda>0$, and define its probability generating function by
\[
 g_\lambda(s)=\exp\{\lambda(s-1)\},\qquad 0\leq s\leq1.
\]
We write $g_\lambda^{\circ r}$ for the $r$-fold composition of $g_\lambda$, with $g_\lambda^{\circ0}(s)=s$. An order-$r$ MIPP, denoted by $V^{(r)}$, has a compound-Poisson representation
\begin{equation}\label{eq:mipp-cp-rep}
 V_t^{(r)}\ \overset{d}=\ \sum_{j=1}^{N_t} Z_j^{(r-1)},\qquad r\geq1,
\end{equation}
where the marks $Z_j^{(r-1)}$ are i.i.d., independent of $N$, and have the distribution of $V_1^{(r-1)}$. For $r=1$, the marks are identically equal to one. The probability generating function and mean are therefore
\begin{equation}\label{eq:mipp-pgf-mean}
 \E[s^{V_1^{(r)}}]=g_\lambda^{\circ r}(s),\qquad
 \E[V_t^{(r)}]=\lambda^r t.
\end{equation}
For $r\geq1$, set
\[
 q_{r-1}=\Pbb(V_1^{(r-1)}>0)=1-g_\lambda^{\circ(r-1)}(0),
\]
where $V_1^{(0)}=1$ and $q_0=1$. Zero marks do not change the count, so positive jumps of $V^{(r)}$ arrive according to a Poisson process of rate $\lambda q_{r-1}$. Conditional on a positive jump, its size has distribution
\begin{equation}\label{eq:mipp-positive-mark}
 \Pbb(Z_j^{(r-1)}=k\mid Z_j^{(r-1)}>0)
 =\frac{\Pbb(V_1^{(r-1)}=k)}{q_{r-1}},\qquad k\geq1.
\end{equation}
We use tildes for an independent MIPP family with base rate $\lambdatilde>0$; its corresponding quantities are $\gtilde$, $\qtilde_j$, and $\Vtilde^{(j)}$.

Let $\{\xi_i\}_{i\geq1}$ and $\{\xitilde_i\}_{i\geq1}$ be independent sequences of strictly positive i.i.d. random variables with distribution functions $F$ and $\Ftilde$, respectively, and finite means
\[
 \mu=\E[\xi_1],\qquad \mutilde=\E[\xitilde_1].
\]
Assume that these sequences are independent of the two MIPP families. For fixed positive integers $m,n$, define the surplus process
\begin{equation}\label{eq:surplus-mn}
 X_t^{(m,n)}(x)
 =x+\sum_{i=1}^{\Vtilde_t^{(m)}}\xitilde_i
      -\sum_{i=1}^{V_t^{(n)}}\xi_i,
 \qquad t\geq0,\quad x\geq0.
\end{equation}
The first sum is aggregate premium income and the second is aggregate claim outgo. Since the expected counts are $\lambdatilde^m t$ and $\lambda^n t$, the expected net surplus drift per unit time is
\[
 \lambdatilde^m\mutilde-\lambda^n\mu.
\]
We impose the net profit condition
\begin{equation}\label{eq:safety-loading-mn}
 \lambdatilde^m\mutilde>\lambda^n\mu.
\end{equation}
Let
\[
 \tau_{m,n}(x)=\inf\{t\geq0:X_t^{(m,n)}(x)<0\},\qquad
 \Lambda_{m,n}(x)=\Pbb(\tau_{m,n}(x)=\infty),
\]
with $\inf\varnothing=\infty$. For a probability distribution $H$, $H^{*k}$ denotes its $k$-fold convolution.

\begin{theorem}[Survival equation and adjustment-coefficient bound]
\label{thm:survival-mn}
Under the assumptions above, the ultimate survival probability satisfies, for every $x\geq0$,
\begin{align}
 (\lambdatilde\qtilde_{m-1}+\lambda q_{n-1})\Lambda_{m,n}(x)
 ={}&\lambdatilde\sum_{k=1}^{\infty}
      \Pbb(\Vtilde_1^{(m-1)}=k)
      \int_{[0,\infty)}\Lambda_{m,n}(x+y)\,d\Ftilde^{*k}(y)\notag\\
 &+\lambda\sum_{k=1}^{\infty}
      \Pbb(V_1^{(n-1)}=k)
      \int_{[0,x]}\Lambda_{m,n}(x-y)\,dF^{*k}(y).
 \label{eq:survival-integral-mn}
\end{align}
Suppose, in addition, that there exists $R>0$ such that $\E[e^{R\xi_1}]<\infty$ and
\begin{equation}\label{eq:adjustment-mn}
 \lambdatilde\left[\gtilde^{\circ(m-1)}
       \left(\E[e^{-R\xitilde_1}]\right)-1\right]
 +\lambda\left[g_\lambda^{\circ(n-1)}
       \left(\E[e^{R\xi_1}]\right)-1\right]=0.
\end{equation}
Then
\begin{equation}\label{eq:martingale-mn}
 M_t=\exp\left\{-R\left(
       \sum_{i=1}^{\Vtilde_t^{(m)}}\xitilde_i
       -\sum_{i=1}^{V_t^{(n)}}\xi_i\right)\right\},\qquad t\geq0,
\end{equation}
is a nonnegative martingale, and
\begin{equation}\label{eq:survival-lower-bound-mn}
 \boxed{\Lambda_{m,n}(x)\geq1-e^{-Rx},\qquad x\geq0.}
\end{equation}
\end{theorem}

\begin{proof}
The positive premium-cluster arrivals form a Poisson process of rate $\lambdatilde\qtilde_{m-1}$, and the positive claim-cluster arrivals form an independent Poisson process of rate $\lambda q_{n-1}$. At a positive premium-cluster arrival, the aggregate premium amount has distribution
\[
 \Gtilde_m(dy)=\frac{1}{\qtilde_{m-1}}
 \sum_{k=1}^{\infty}\Pbb(\Vtilde_1^{(m-1)}=k)\,d\Ftilde^{*k}(y),
\]
while the aggregate claim amount at a positive claim-cluster arrival has distribution
\[
 G_n(dy)=\frac{1}{q_{n-1}}
 \sum_{k=1}^{\infty}\Pbb(V_1^{(n-1)}=k)\,dF^{*k}(y).
\]
The first nonzero cluster arrival is exponential with rate $\lambdatilde\qtilde_{m-1}+\lambda q_{n-1}$. Conditioning on its type and size gives
\begin{align*}
 &(\lambdatilde\qtilde_{m-1}+\lambda q_{n-1})\Lambda_{m,n}(x)\\
 &\quad=\lambdatilde\qtilde_{m-1}\int_{[0,\infty)}\Lambda_{m,n}(x+y)\,\Gtilde_m(dy)
 +\lambda q_{n-1}\int_{[0,x]}\Lambda_{m,n}(x-y)\,G_n(dy).
\end{align*}
Substitution of the two cluster-size mixtures proves \eqref{eq:survival-integral-mn}.

For the martingale assertion, let
\[
 \widetilde W_m=\sum_{j=1}^{\Vtilde_1^{(m-1)}}\xitilde_j,
 \qquad W_n=\sum_{j=1}^{V_1^{(n-1)}}\xi_j,
\]
with an empty sum interpreted as zero. The generating-function identities give
\[
 \E[e^{-R\widetilde W_m}]
 =\gtilde^{\circ(m-1)}(\E[e^{-R\xitilde_1}]),\qquad
 \E[e^{RW_n}]
 =g_\lambda^{\circ(n-1)}(\E[e^{R\xi_1}]).
\]
The compound-Poisson exponential formula and independence yield
\[
 \E[M_t]=\exp\left\{t\left(\lambdatilde(\E[e^{-R\widetilde W_m}]-1)
 +\lambda(\E[e^{RW_n}]-1)\right)\right\}=1.
\]
Independent increments then show that $M$ is a martingale. Let $T>0$. Optional sampling at the bounded stopping time $\tau_{m,n}(x)\wedge T$ gives
\[
 1=\E[M_{\tau_{m,n}(x)\wedge T}].
\]
On $\{\tau_{m,n}(x)\leq T\}$, the aggregate net payments are less than $-x$, so $M_{\tau_{m,n}(x)}>e^{Rx}$. Nonnegativity therefore implies
\[
 1\geq e^{Rx}\Pbb(\tau_{m,n}(x)\leq T).
\]
Letting $T\to\infty$ proves the bound.
\end{proof}

\begin{remark}
The positive net profit condition implies that the derivative at zero of the adjustment function in \eqref{eq:adjustment-mn} is negative:
\[
 \psi'_{m,n}(0+)=-\lambdatilde^m\mutilde+\lambda^n\mu<0.
\]
It does not by itself guarantee that a positive root exists. The exponential bound is asserted only when a positive solution of the adjustment equation exists and the relevant exponential moment is finite at that solution. The strong law of large numbers also implies $\Lambda_{m,n}(x)\to1$ as $x\to\infty$ under \eqref{eq:safety-loading-mn}.
\end{remark}

\subsection{Premium arrivals follow an ordinary Poisson process}
\label{subsec:poisson-premium}

We first specialize to ordinary Poisson premium arrivals, while claim counts continue to be described by an order-$n$ MIPP. Let $\Ntilde$ be a Poisson process with rate $\lambdatilde$, independent of $V^{(n)}$, and set
\begin{equation}\label{eq:surplus-n}
 X_t^{(n)}(x)=x+\sum_{i=1}^{\Ntilde_t}\xitilde_i
 -\sum_{i=1}^{V_t^{(n)}}\xi_i,
 \qquad t\geq0,\quad x\geq0.
\end{equation}
The net profit condition is
\begin{equation}\label{eq:safety-loading-n}
 \lambdatilde\mutilde>\lambda^n\mu.
\end{equation}
Write $\Lambda_n(x)$ for the ultimate survival probability of this model.

\begin{corollary}[Survival equation and adjustment-coefficient bound]
\label{cor:survival-n}
For all $x\geq0$,
\begin{align}
 (\lambda q_{n-1}+\lambdatilde)\Lambda_n(x)
 ={}&\lambda\sum_{k=1}^{\infty}\Pbb(V_1^{(n-1)}=k)
 \int_{[0,x]}\Lambda_n(x-y)\,dF^{*k}(y)\notag\\
 &+\lambdatilde\int_{[0,\infty)}\Lambda_n(x+y)\,d\Ftilde(y).
 \label{eq:survival-n}
\end{align}
If $R_n>0$ satisfies
\begin{equation}\label{eq:adjustment-n}
 \lambdatilde(\E[e^{-R_n\xitilde_1}]-1)
 +\lambda\left[g_\lambda^{\circ(n-1)}(\E[e^{R_n\xi_1}])-1\right]=0,
\end{equation}
with the required exponential moment finite, then the corresponding exponential process is a martingale and
\begin{equation}\label{eq:bound-n}
 \Lambda_n(x)\geq1-e^{-R_nx},\qquad x\geq0.
\end{equation}
\end{corollary}
\begin{proof}
This is the case $m=1$ of Theorem~\ref{thm:survival-mn}; in this case $\Vtilde^{(1)}=\Ntilde$ and $\qtilde_0=1$.
\end{proof}

We next obtain a transform formula when the individual premium amount is exponentially distributed.

\begin{theorem}[Laplace transform for exponential premium amounts]
\label{thm:laplace-survival-n}
Suppose $\xitilde_1$ is exponentially distributed with mean $\mutilde$, while $\xi_1$ has finite mean $\mu$. Define
\[
 \cL_{\Lambda_n}(\theta)=\int_0^\infty e^{-\theta x}\Lambda_n(x)\,dx,
 \qquad \widehat F(\theta)=\E[e^{-\theta\xi_1}],\quad \theta>0.
\]
Under \eqref{eq:safety-loading-n},
\begin{equation}\label{eq:laplace-survival-n}
 \boxed{\displaystyle
 \cL_{\Lambda_n}(\theta)=
 \frac{\lambdatilde\mutilde-\lambda^n\mu}
 {\lambda(1-\mutilde\theta)
 [g_\lambda^{\circ(n-1)}(\widehat F(\theta))-1]
 +\lambdatilde\mutilde\theta}.}
\end{equation}
At $\theta=1/\mutilde$, the right-hand side is understood by continuity.
\end{theorem}

\begin{proof}
The density of $\xitilde_1$ is $\mutilde^{-1}e^{-y/\mutilde}$ on $(0,\infty)$. Taking Laplace transforms in \eqref{eq:survival-n} and applying Tonelli's theorem gives
\[
 (\lambda q_{n-1}+\lambdatilde)\cL_{\Lambda_n}(\theta)
 =\lambda\cL_{\Lambda_n}(\theta)
 \sum_{k=1}^{\infty}\Pbb(V_1^{(n-1)}=k)\widehat F(\theta)^k
 +\frac{\lambdatilde}{\mutilde}I_\theta,
\]
where
\[
 I_\theta=\int_0^\infty\int_0^\infty
 \Lambda_n(x+y)e^{-\theta x-y/\mutilde}\,dy\,dx.
\]
For $\theta\ne1/\mutilde$, a change of variables gives
\[
 \frac{1}{\mutilde}I_\theta
 =\frac{\cL_{\Lambda_n}(\theta)-\cL_{\Lambda_n}(1/\mutilde)}
 {1-\mutilde\theta}.
\]
Let $p_0=\Pbb(V_1^{(n-1)}=0)=1-q_{n-1}$. The probability generating function gives
\[
 \sum_{k=1}^{\infty}\Pbb(V_1^{(n-1)}=k)\widehat F(\theta)^k
 =g_\lambda^{\circ(n-1)}(\widehat F(\theta))-p_0.
\]
After rearrangement,
\begin{equation}\label{eq:laplace-intermediate-n}
 \lambdatilde\cL_{\Lambda_n}(1/\mutilde)
 =\left\{\lambda(1-\mutilde\theta)
 [g_\lambda^{\circ(n-1)}(\widehat F(\theta))-1]
 +\lambdatilde\mutilde\theta\right\}\cL_{\Lambda_n}(\theta).
\end{equation}
By the net profit condition and the strong law, $\Lambda_n(x)\to1$ as $x\to\infty$, and hence
\[
 \lim_{\theta\downarrow0}\theta\cL_{\Lambda_n}(\theta)=1.
\]
Further,
\[
 g_\lambda^{\circ(n-1)}(1)=1,\quad
 (g_\lambda^{\circ(n-1)})'(1)=\lambda^{n-1},\quad
 \widehat F(\theta)=1-\mu\theta+o(\theta).
\]
Thus $g_\lambda^{\circ(n-1)}(\widehat F(\theta))-1
=-\lambda^{n-1}\mu\theta+o(\theta)$. Letting $\theta\downarrow0$ in \eqref{eq:laplace-intermediate-n} yields
\[
 \lambdatilde\cL_{\Lambda_n}(1/\mutilde)=\lambdatilde\mutilde-\lambda^n\mu.
\]
Substituting this identity into \eqref{eq:laplace-intermediate-n} proves the theorem. Continuity gives the value at $\theta=1/\mutilde$.
\end{proof}

\begin{example}[Exponential premiums and claims]
\label{ex:poisson-exponential}
When $n=1$ and both individual premium and claim amounts are exponential with means $\mutilde$ and $\mu$, respectively, $g_\lambda^{\circ0}(s)=s$ and $\widehat F(\theta)=(1+\mu\theta)^{-1}$. Formula~\eqref{eq:laplace-survival-n} reduces to
\[
 \cL_{\Lambda_1}(\theta)=
 \frac{(\lambdatilde\mutilde-\lambda\mu)(1+\mu\theta)}
 {\theta\{\lambdatilde\mutilde-\lambda\mu
 +(\lambda+\lambdatilde)\mu\mutilde\theta\}}.
\]
Inverting this transform gives
\begin{equation}\label{eq:poisson-exp-survival}
 \Lambda_1(x)=1-
 \frac{\lambda(\mu+\mutilde)}{(\lambda+\lambdatilde)\mutilde}
 \exp\left\{-\frac{\lambdatilde\mutilde-\lambda\mu}
 {(\lambda+\lambdatilde)\mu\mutilde}x\right\},\qquad x\geq0.
\end{equation}
The exponent is negative under the net profit condition.
\end{example}

\begin{theorem}[Limit as the claim MIPP order increases]
\label{thm:claim-order-limit}
Assume the conditions of Theorem~\ref{thm:laplace-survival-n} and, in addition, $\lambdatilde\mutilde>\mu$. Then, for each fixed $x>0$,
\begin{equation}\label{eq:claim-order-limit}
 \lim_{n\to\infty}\Lambda_n(x)=
 \begin{cases}
 0,&\lambda>1,\\[1mm]
 1-\displaystyle\frac{\mu}{\lambdatilde\mutilde},&\lambda=1,\\[2mm]
 1,&0<\lambda<1.
 \end{cases}
\end{equation}
\end{theorem}
\begin{proof}
If $\lambda>1$, then $\lambda^n\mu\to\infty$, so the net drift is negative for all sufficiently large $n$. The strong law implies $X_t^{(n)}(x)\to-\infty$ almost surely, and therefore $\Lambda_n(x)=0$ eventually.

If $\lambda=1$, then $g_1(s)=e^{s-1}$ has 1 as its unique fixed point in $[0,1]$ and $g_1(s)>s$ for $s<1$. Hence $g_1^{\circ(n-1)}(s)\to1$ for $s\in[0,1)$. Formula~\eqref{eq:laplace-survival-n} then gives
\[
 \cL_{\Lambda_n}(\theta)\longrightarrow
 \frac{\lambdatilde\mutilde-\mu}{\lambdatilde\mutilde\theta}.
\]
For $0<\lambda<1$, the same fixed-point argument applies to $g_\lambda$, while $\lambda^n\mu\to0$, and therefore
\[
 \cL_{\Lambda_n}(\theta)\longrightarrow\frac1\theta.
\]
To pass from transform convergence to pointwise convergence, note that $\Lambda_n$ is nondecreasing and bounded between 0 and 1. Helly's selection theorem implies that every subsequence has a further subsequence converging at all continuity points of a nondecreasing limit $\Lambda_*$. Dominated convergence identifies its Laplace transform as $c/\theta$, where $c=1-\mu/(\lambdatilde\mutilde)$ in the critical case and $c=1$ when $\lambda<1$. Uniqueness of Laplace transforms yields $\Lambda_*(x)=c$ almost everywhere. Monotonicity and a squeeze by continuity points on either side of a fixed $x>0$ then give convergence at that $x$. Since every subsequence has a further subsequence with the same limit, the whole sequence converges pointwise.
\end{proof}

\subsection{Premium counts follow a compound-Poisson process}
\label{subsec:compound-poisson-premium}

We next consider a genuinely clustered premium-count process. For clarity, we take the premium count to be an order-two MIPP,
\[
 \Vtilde_t^{(2)}=\sum_{j=1}^{\Ntilde_t}Z_j,
 \qquad \Ntilde_t\sim\operatorname{Poisson}(\lambdatilde t),
 \qquad Z_j\sim\operatorname{Poisson}(\lambdatilde),
\]
where the outer process and marks are independent. Thus $\E[\Vtilde_t^{(2)}]=\lambdatilde^2t$. The surplus process is
\begin{equation}\label{eq:surplus-compound-premium}
 X_t^{(2,n)}(x)=x+\sum_{i=1}^{\Vtilde_t^{(2)}}\xitilde_i
 -\sum_{i=1}^{V_t^{(n)}}\xi_i.
\end{equation}
The net profit condition is
\begin{equation}\label{eq:safety-compound-premium}
 \lambdatilde^2\mutilde>\lambda^n\mu.
\end{equation}
Positive premium clusters arrive at rate
\[
 a_2=\lambdatilde(1-e^{-\lambdatilde}),
\]
and the number of individual premium payments within a positive cluster has a Poisson$(\lambdatilde)$ distribution conditioned to be positive. The aggregate premium-cluster distribution is therefore
\begin{equation}\label{eq:premium-cluster-law}
 H_2(dy)=\frac{1}{1-e^{-\lambdatilde}}
 \sum_{k=1}^{\infty}e^{-\lambdatilde}
 \frac{\lambdatilde^k}{k!}\,\Ftilde^{*k}(dy).
\end{equation}
For claims, define
\[
 G_n(dy)=\frac1{q_{n-1}}
 \sum_{k=1}^{\infty}\Pbb(V_1^{(n-1)}=k)F^{*k}(dy).
\]

\begin{corollary}[Survival equation and ruin bound for clustered premiums]
\label{cor:compound-premium}
The ultimate survival probability $\Lambda_{2,n}(x)$ satisfies
\begin{align}
 (a_2+\lambda q_{n-1})\Lambda_{2,n}(x)
 ={}&a_2\int_{[0,\infty)}\Lambda_{2,n}(x+y)H_2(dy)\notag\\
 &+\lambda\sum_{k=1}^{\infty}\Pbb(V_1^{(n-1)}=k)
 \int_{[0,x]}\Lambda_{2,n}(x-y)F^{*k}(dy).
 \label{eq:compound-premium-renewal}
\end{align}
If there is $R>0$ with $\E[e^{R\xi_1}]<\infty$ such that
\begin{equation}\label{eq:compound-premium-adjustment}
 \lambdatilde\left[g_{\lambdatilde}(\E[e^{-R\xitilde_1}])-1\right]
 +\lambda\left[g_\lambda^{\circ(n-1)}(\E[e^{R\xi_1}])-1\right]=0,
\end{equation}
then the associated exponential process is a nonnegative martingale and
\begin{equation}\label{eq:compound-premium-bound}
 \Lambda_{2,n}(x)\geq1-e^{-Rx},\qquad x\geq0.
\end{equation}
\end{corollary}
\begin{proof}
The first identity follows by conditioning on the first positive premium or claim cluster arrival, whose total rate is $a_2+\lambda q_{n-1}$. The cluster laws are $H_2$ and $G_n$, respectively, giving \eqref{eq:compound-premium-renewal}. For the martingale, the compound-Poisson exponential formula gives the premium contribution
\[
 \E\left[e^{-R\sum_{i=1}^{\Vtilde_t^{(2)}}\xitilde_i}\right]
 =\exp\left\{t\lambdatilde\left[g_{\lambdatilde}(\E[e^{-R\xitilde_1}])-1\right]\right\},
\]
and the claim contribution
\[
 \E\left[e^{R\sum_{i=1}^{V_t^{(n)}}\xi_i}\right]
 =\exp\left\{t\lambda\left[g_\lambda^{\circ(n-1)}(\E[e^{R\xi_1}])-1\right]\right\}.
\]
Equation~\eqref{eq:compound-premium-adjustment} makes their product equal to one. Independent increments yield the martingale property, and the optional-stopping argument used in Theorem~\ref{thm:survival-mn} proves the bound.
\end{proof}

\begin{remark}[Individual and aggregate exponential premium sizes]
Even if the individual premium amount $\xitilde_1$ is exponentially distributed, the aggregate positive-cluster amount in \eqref{eq:premium-cluster-law} is generally a mixture of Erlang distributions, not an exponential distribution. Therefore, the simple transform formula below requires an explicit exponential assumption on the aggregate premium-cluster amount; it does not follow merely from exponential individual premium amounts.
\end{remark}

\begin{theorem}[Transform under an exponential aggregate premium-cluster law]
\label{thm:compound-premium-laplace}
Suppose the aggregate premium amount at each positive premium-cluster arrival is exponential with mean $\bar\mu_2$, independently across clusters and of the claim process. If $a_2\bar\mu_2>\lambda^n\mu$, then
\begin{equation}\label{eq:compound-premium-laplace}
 \boxed{\displaystyle
 \cL_{\Lambda_{2,n}}(\theta)=
 \frac{a_2\bar\mu_2-\lambda^n\mu}
 {\lambda(1-\bar\mu_2\theta)
 [g_\lambda^{\circ(n-1)}(\hatF(\theta))-1]
 +a_2\bar\mu_2\theta},\qquad\theta>0,}
\end{equation}
where $\hatF(\theta)=\E[e^{-\theta\xi_1}]$. At $\theta=1/\bar\mu_2$, the expression is interpreted by continuity.
\end{theorem}
\begin{proof}
The proof is the same Laplace-transform calculation as in Theorem~\ref{thm:laplace-survival-n}, with the positive premium-cluster rate $a_2$ and exponential cluster mean $\bar\mu_2$ replacing $\lambdatilde$ and $\mutilde$. In particular, the transformed premium integral is
\[
 \frac{\cL_{\Lambda_{2,n}}(\theta)-\cL_{\Lambda_{2,n}}(1/\bar\mu_2)}
 {1-\bar\mu_2\theta},
\]
and the small-$\theta$ expansion of the claim generating function is
\[
 g_\lambda^{\circ(n-1)}(\hatF(\theta))-1
 =-\lambda^{n-1}\mu\theta+o(\theta).
\]
The net profit condition gives $\theta\cL_{\Lambda_{2,n}}(\theta)\to1$, which identifies the transform value at $1/\bar\mu_2$ as $a_2\bar\mu_2-\lambda^n\mu$. Substitution yields \eqref{eq:compound-premium-laplace}.
\end{proof}

\begin{remark}[Mean-matched cluster approximation]
For the original order-two MIPP with individual premium payments, the mean aggregate amount conditional on a positive cluster is
\[
 \E[W_2]=\frac{\lambdatilde\mutilde}{1-e^{-\lambdatilde}}.
\]
A mean-matched exponential aggregate-cluster model is obtained by choosing $\bar\mu_2=\lambdatilde\mutilde/(1-e^{-\lambdatilde})$. This matches the expected premium inflow, since $a_2\bar\mu_2=\lambdatilde^2\mutilde$, but does not make the full cluster-size distributions identical.
\end{remark}

\section{Models with constant and stochastic premium income}
\label{sec:constant-plus-stochastic}

We now incorporate a deterministic premium rate in addition to stochastic premium payments. Consider
\begin{equation}\label{eq:constant-stochastic-surplus}
 X_t(x)=x+ct+\sum_{i=1}^{\Ntilde_t}\xitilde_i
 -\sum_{i=1}^{N_t}\xi_i,
 \qquad c>0,\quad t\geq0,\quad x\geq0,
\end{equation}
where $N$ and $\Ntilde$ are independent Poisson processes with rates $\lambda$ and $\lambdatilde$. The claim and premium sizes are independent, strictly positive, and have distributions $F$ and $\Ftilde$, with finite means $\mu$ and $\mutilde$. The net profit condition is
\begin{equation}\label{eq:constant-stochastic-safety}
 c+\lambdatilde\mutilde>\lambda\mu.
\end{equation}
Let $\tau(x)=\inf\{t\geq0:X_t(x)<0\}$ and $\Lambda(x)=\Pbb(\tau(x)=\infty)$.

\begin{theorem}[Integro-differential equation and adjustment-coefficient bound]
\label{thm:constant-stochastic}
Assume that $F$ and $\Ftilde$ are continuous. Then, for $x>0$, the survival probability satisfies
\begin{equation}\label{eq:survival-ide}
 c\Lambda'(x)=(\lambda+\lambdatilde)\Lambda(x)
 -\lambda\int_{[0,x]}\Lambda(x-y)\,dF(y)
 -\lambdatilde\int_{[0,\infty)}\Lambda(x+y)\,d\Ftilde(y).
\end{equation}
If there exists $R>0$ for which $\E[e^{R\xi}]<\infty$ and
\begin{equation}\label{eq:constant-stochastic-adjustment}
 -Rc+\lambdatilde(\E[e^{-R\xitilde}]-1)
 +\lambda(\E[e^{R\xi}]-1)=0,
\end{equation}
then
\[
 M_t=\exp\left\{-R\left(ct+
 \sum_{i=1}^{\Ntilde_t}\xitilde_i-\sum_{i=1}^{N_t}\xi_i\right)\right\}
\]
is a nonnegative martingale and
\begin{equation}\label{eq:constant-stochastic-bound}
 \Lambda(x)\geq1-e^{-Rx},\qquad x\geq0.
\end{equation}
\end{theorem}

\begin{proof}
For a small time interval, the probability of no premium or claim arrival is $1-(\lambda+\lambdatilde)t+o(t)$; the probability of one claim arrival is $\lambda t+o(t)$; and the probability of one premium arrival is $\lambdatilde t+o(t)$. The probability of two or more arrivals is $o(t)$. Conditioning on these events and using the deterministic increase $ct$ gives the backward equation
\[
\begin{split}
\Lambda(x)={}&[1-(\lambda+\lambdatilde)t]\Lambda(x+ct)\\
&+\lambda t\int_{[0,x+ct]}\Lambda(x+ct-y)\,dF(y)\\
&+\lambdatilde t\int_{[0,\infty)}\Lambda(x+ct+y)\,d\Ftilde(y)+o(t).
\end{split}
\]
Dividing by $t$ and letting $t\downarrow0$ yields \eqref{eq:survival-ide}. For the martingale, the compound-Poisson exponential formula gives
\[
 \E[M_t]=\exp\{t[-Rc+\lambdatilde(\E[e^{-R\xitilde}]-1)
 +\lambda(\E[e^{R\xi}]-1)]\}=1.
\]
Independent increments imply that $M$ is a martingale. Optional sampling at $\tau(x)\wedge T$ and the fact that $X_{\tau(x)}(x)<0$ on ruin give
\[
 1=\E[M_{\tau(x)\wedge T}]\geq e^{Rx}\Pbb(\tau(x)\leq T).
\]
Letting $T\to\infty$ proves \eqref{eq:constant-stochastic-bound}.
\end{proof}

\begin{remark}[Ruin-probability equation]
Let $\Pi(x)=1-\Lambda(x)$ denote the ultimate ruin probability. Since the premium and claim amounts are strictly positive, $F(0)=\Ftilde(0)=0$. Equation~\eqref{eq:survival-ide} is equivalently
\begin{equation}\label{eq:ruin-ide-corrected}
 -c\Pi'(x)=\lambda[1-F(x)]-(\lambda+\lambdatilde)\Pi(x)
 +\lambda\int_{[0,x]}\Pi(x-y)\,dF(y)
 +\lambdatilde\int_{[0,\infty)}\Pi(x+y)\,d\Ftilde(y).
\end{equation}
The term $-(\lambda+\lambdatilde)\Pi(x)$ is required when passing from the survival equation to the ruin-probability equation.
\end{remark}

\subsection{Exponential premium amounts}
\label{subsec:exponential-premiums}

Suppose now that individual premium amounts are exponential with mean $\mutilde$, while claim amounts retain a general distribution $F$ with finite mean $\mu$. Define
\[
 \cL_\Lambda(\theta)=\int_0^\infty e^{-\theta x}\Lambda(x)\,dx,
 \qquad \hatF(\theta)=\int_{[0,\infty)}e^{-\theta y}\,dF(y).
\]

\begin{theorem}[Laplace transform of the survival probability]
\label{thm:laplace-constant-stochastic}
Under \eqref{eq:constant-stochastic-safety},
\begin{equation}\label{eq:laplace-constant-stochastic}
 \boxed{\displaystyle
 \cL_\Lambda(\theta)=
 \frac{c+\lambdatilde\mutilde-\lambda\mu
       -c\mutilde\theta\Lambda(0+)}
 {(c\theta-\lambda+\lambda\hatF(\theta))(1-\mutilde\theta)
       +\lambdatilde\mutilde\theta},\qquad\theta>0.}
\end{equation}
At $\theta=1/\mutilde$, the expression is interpreted by continuity.
\end{theorem}

\begin{proof}
Multiply \eqref{eq:survival-ide} by $e^{-\theta x}$ and integrate over $x>0$. Integration by parts gives
\[
 \int_0^\infty e^{-\theta x}\Lambda'(x)\,dx
 =-\Lambda(0+)+\theta\cL_\Lambda(\theta).
\]
Tonelli's theorem gives the transform of the claim convolution as $\cL_\Lambda(\theta)\hatF(\theta)$. For the premium term, the exponential density yields
\[
 \int_0^\infty e^{-\theta x}\int_0^\infty
 \Lambda(x+y)\frac{e^{-y/\mutilde}}{\mutilde}\,dy\,dx
 =\frac{\cL_\Lambda(\theta)-\cL_\Lambda(1/\mutilde)}
 {1-\mutilde\theta}.
\]
Substitution and rearrangement express $\cL_\Lambda(\theta)$ in terms of $\cL_\Lambda(1/\mutilde)$. Since positive drift implies $\Lambda(x)\to1$, we have $\theta\cL_\Lambda(\theta)\to1$ as $\theta\downarrow0$. Also,
\[
 \lim_{\theta\downarrow0}\frac{\hatF(\theta)-1}{\theta}=-\mu.
\]
Taking the limit in the rearranged identity gives
\[
 \lambdatilde\cL_\Lambda(1/\mutilde)
 =c[1-\Lambda(0+)]+\lambdatilde\mutilde-\lambda\mu.
\]
Substituting this value yields \eqref{eq:laplace-constant-stochastic}.
\end{proof}

\subsection{Exponential premium and claim amounts}
\label{subsec:both-exponential}

Assume in addition that individual claim amounts are exponential with mean $\mu$. Then $\hatF(\theta)=(1+\mu\theta)^{-1}$, and the denominator in \eqref{eq:laplace-constant-stochastic} is a cubic polynomial with a factor $\theta$. Let
\[
 C=c+\lambdatilde\mutilde-\lambda\mu,
 \qquad A=(\lambda+\lambdatilde)\mutilde\mu-c(\mutilde-\mu).
\]
The two nonzero roots are the roots of
\begin{equation}\label{eq:theta-roots}
 \theta^2-\frac{A}{c\mutilde\mu}\theta
 -\frac{C}{c\mutilde\mu}=0,
\end{equation}
namely
\begin{equation}\label{eq:theta12}
 \theta_1=\frac{A+\sqrt{A^2+4c\mutilde\mu C}}{2c\mutilde\mu}>0,
 \qquad
 \theta_2=\frac{A-\sqrt{A^2+4c\mutilde\mu C}}{2c\mutilde\mu}<0.
\end{equation}

\begin{theorem}[Initial survival probability and explicit survival function]
\label{thm:explicit-constant-stochastic}
Under the assumptions above,
\begin{equation}\label{eq:lambda-zero}
 \boxed{\displaystyle
 \Lambda(0+)=\frac{c+\lambdatilde\mutilde-\lambda\mu}
 {c\mutilde\theta_1}=-\mu\theta_2\in(0,1).}
\end{equation}
Moreover,
\begin{equation}\label{eq:explicit-constant-stochastic}
 \boxed{\displaystyle
 \Lambda(x)=1-\bigl(1-\Lambda(0+)\bigr)e^{\theta_2x},\qquad x\geq0.}
\end{equation}
\end{theorem}
\begin{proof}
After substituting $\hatF(\theta)=(1+\mu\theta)^{-1}$ into \eqref{eq:laplace-constant-stochastic} and multiplying numerator and denominator by $1+\mu\theta$, the denominator factors as $\theta$ times a quadratic with roots $\theta_1$ and $\theta_2$. Since $\cL_\Lambda$ is finite and analytic for $\theta>0$, the positive root $\theta_1$ must be cancelled by the numerator. This gives
\[
 c\mutilde\theta_1\Lambda(0+)=C.
\]
The product of the roots in \eqref{eq:theta-roots} is $\theta_1\theta_2=-C/(c\mutilde\mu)$, which yields $\Lambda(0+)=-\mu\theta_2$. Since $C>0$, we have $\theta_2<0$. Furthermore, substituting $-1/\mu$ into the quadratic in \eqref{eq:theta-roots} gives a positive value, while its value at zero is negative; hence $-1/\mu<\theta_2<0$ and $0<\Lambda(0+)<1$. Cancelling the factor $\theta-\theta_1$ leaves
\[
 \cL_\Lambda(\theta)=\frac1\theta
 -\frac{1-\Lambda(0+)}{\theta-\theta_2}.
\]
Laplace inversion proves \eqref{eq:explicit-constant-stochastic}.
\end{proof}

\subsection{MIPP version with a constant premium rate}
\label{subsec:mipp-constant}

The preceding model extends naturally when the stochastic premium and claim counts are both MIPP-driven. Suppose the individual premium and claim amounts have continuous distributions, as above. Consider
\begin{equation}\label{eq:mipp-constant-surplus}
 X_t^{(m,n)}(x)=x+ct+
 \sum_{i=1}^{\Vtilde_t^{(m)}}\xitilde_i
 -\sum_{i=1}^{V_t^{(n)}}\xi_i,
 \qquad c>0.
\end{equation}
The net profit condition is
\begin{equation}\label{eq:mipp-constant-safety}
 c+\lambdatilde^m\mutilde>\lambda^n\mu.
\end{equation}
Let
\[
 \tau^c_{m,n}(x)=\inf\{t\geq0:X_t^{(m,n)}(x)<0\},
 \qquad \Lambda^c_{m,n}(x)=\Pbb(\tau^c_{m,n}(x)=\infty)
\]
be the ruin time and ultimate survival probability for this constant-plus-stochastic premium model.

\begin{theorem}[MIPP integro-differential equation and exponential bound]
\label{thm:mipp-constant}
For $x>0$, the survival function satisfies
\begin{align}
 c(\Lambda^c_{m,n})'(x)
 ={}&(\lambdatilde\qtilde_{m-1}+\lambda q_{n-1})\Lambda^c_{m,n}(x)\notag\\
 &-\lambdatilde\sum_{k=1}^{\infty}\Pbb(\Vtilde_1^{(m-1)}=k)
 \int_{[0,\infty)}\Lambda^c_{m,n}(x+y)\,d\Ftilde^{*k}(y)\notag\\
 &-\lambda\sum_{k=1}^{\infty}\Pbb(V_1^{(n-1)}=k)
 \int_{[0,x]}\Lambda^c_{m,n}(x-y)\,dF^{*k}(y).
 \label{eq:mipp-constant-ide}
\end{align}
If $R>0$ solves
\begin{equation}\label{eq:mipp-constant-adjustment}
 -Rc+\lambdatilde\left[\gtilde^{\circ(m-1)}(\E[e^{-R\xitilde_1}])-1\right]
 +\lambda\left[g_\lambda^{\circ(n-1)}(\E[e^{R\xi_1}])-1\right]=0,
\end{equation}
with the required exponential moment finite, then the corresponding exponential process is a nonnegative martingale and
\begin{equation}\label{eq:mipp-constant-bound}
 \Lambda^c_{m,n}(x)\geq1-e^{-Rx},\qquad x\geq0.
\end{equation}
\end{theorem}
\begin{proof}
The positive premium and claim cluster rates are $\lambdatilde\qtilde_{m-1}$ and $\lambda q_{n-1}$. Conditioning over an infinitesimal time interval, during which the surplus increases deterministically at rate $c$, gives \eqref{eq:mipp-constant-ide}. The exponential formula for the marked compound-Poisson representations shows that the process
\[
 \exp\left\{-R\left(ct+
 \sum_{i=1}^{\Vtilde_t^{(m)}}\xitilde_i-
 \sum_{i=1}^{V_t^{(n)}}\xi_i\right)\right\}
\]
has expectation one and independent multiplicative increments under \eqref{eq:mipp-constant-adjustment}; it is therefore a martingale. Optional sampling at the ruin time truncated at $T$ gives \eqref{eq:mipp-constant-bound}, as in Theorem~\ref{thm:constant-stochastic}.
\end{proof}

\section{Finite-support claim clusters: explicit illustrations}
\label{sec:examples}

We now consider three model classes in which the number of individual claims generated by one claim event has a finite-support distribution. These examples allow explicit transforms and, under exponential individual claim sizes, rational Laplace transforms that can be inverted by partial fractions.

\subsection{Type I: deterministic premium income and clustered claims}
\label{subsec:type-one}

Let claims arrive according to a Poisson process with rate $\lambda$. At each claim arrival, a random number $V$ of individual claims is generated, where
\[
 \Pbb(V=i)=p_i,\quad i=1,\ldots,n,\qquad \sum_{i=1}^n p_i=1.
\]
The individual claim sizes $\xi_1,\xi_2,\ldots$ are i.i.d. with distribution $F$ and mean $\mu$. The surplus process with constant premium rate $c>0$ is
\begin{equation}\label{eq:type-one-surplus}
 X_t(x)=x+ct-\sum_{j=1}^{N_t}\sum_{k=1}^{V_j}\xi_{j,k},
\end{equation}
where the $V_j$ are i.i.d. copies of $V$, independent of $N$ and the claim sizes. The net profit condition is
\begin{equation}\label{eq:type-one-safety}
 c>\lambda\mu\E[V]=\lambda\mu\sum_{i=1}^n i p_i.
\end{equation}
The aggregate claim-size distribution at one event is $H=\sum_{i=1}^n p_iF^{*i}$, with transform
\[
 \hatH(\theta)=\sum_{i=1}^n p_i\hatF(\theta)^i.
\]

\begin{theorem}[Laplace transform for Type I]
\label{thm:type-one-transform}
The Laplace transform of the ultimate survival probability is
\begin{equation}\label{eq:type-one-transform}
 \boxed{\displaystyle
 \cL_\Lambda(\theta)=
 \frac{c-\lambda\mu\sum_{i=1}^n i p_i}
 {c\theta-\lambda\left(1-\sum_{i=1}^n p_i\hatF(\theta)^i\right)},\qquad\theta>0.}
\end{equation}
\end{theorem}
\begin{proof}
For this classical compound-Poisson risk model, the survival probability satisfies
\[
 c\Lambda'(x)=\lambda\Lambda(x)
 -\lambda\sum_{i=1}^n p_i\int_{[0,x]}\Lambda(x-y)\,dF^{*i}(y).
\]
Taking Laplace transforms yields
\[
 c[-\Lambda(0+)+\theta\cL_\Lambda(\theta)]
 =\lambda\cL_\Lambda(\theta)
 \left[1-\sum_{i=1}^n p_i\hatF(\theta)^i\right].
\]
The standard net-profit identity for this model gives
\[
 \Lambda(0+)=1-\frac{\lambda\mu\sum_{i=1}^n i p_i}{c}.
\]
Substitution proves \eqref{eq:type-one-transform}.
\end{proof}

Assume in the remainder of this subsection that individual claim sizes are exponential with mean $\mu$. Then $\hatF(\theta)=(1+\mu\theta)^{-1}$, and \eqref{eq:type-one-transform} is a rational function. The next result gives a convenient general inversion formula.

\begin{theorem}[Inversion of rational survival transforms]
\label{thm:rational-inversion}
Suppose a rational Laplace transform of a survival function has the partial-fraction representation
\begin{equation}\label{eq:rational-partial-fractions}
 \cL_\Lambda(\theta)=\frac1\theta+
 \sum_{\ell=1}^{L}\sum_{j=1}^{m_\ell}
 \frac{c_{\ell j}}{(\theta-r_\ell)^j},
 \qquad \Re(r_\ell)<0,
\end{equation}
where conjugate poles and their coefficients occur in pairs so that the transform is real-valued for real $\theta>0$. Then
\begin{equation}\label{eq:rational-inversion}
 \boxed{\displaystyle
 \Lambda(x)=1+\sum_{\ell=1}^{L}\sum_{j=1}^{m_\ell}
 c_{\ell j}\frac{x^{j-1}}{(j-1)!}e^{r_\ell x},\qquad x\geq0.}
\end{equation}
Conjugate terms combine to give real expressions involving exponentially damped sine and cosine functions.
\end{theorem}
\begin{proof}
The result follows by linearity of the inverse Laplace transform and the standard identity (see, for example, \cite{PolyaninManzhirov2008})
\[
 \cL^{-1}\left\{\frac{1}{(\theta-r)^j}\right\}(x)
 =\frac{x^{j-1}}{(j-1)!}e^{rx}.
\]
The inverse transform of $1/\theta$ is the constant function one. Summing the terms proves the formula.
\end{proof}

For $n=1$, the model reduces to the classical Cram\'er--Lundberg model. If both claim sizes and cluster sizes are exponential and deterministic one, then
\[
 \cL_\Lambda(\theta)=\frac{c-\lambda\mu}{c\theta-\lambda(1-(1+\mu\theta)^{-1})}
 =\frac1\theta-\frac{\lambda\mu/c}{\theta+(c-\lambda\mu)/(c\mu)},
\]
so that
\begin{equation}\label{eq:type-one-single}
 \Lambda(x)=1-\frac{\lambda\mu}{c}
 \exp\left\{-\frac{c-\lambda\mu}{c\mu}x\right\}.
\end{equation}

\begin{example}[One, two, and three possible claims per event]
\label{ex:type-one-examples}
For $p_1=1$, $\lambda=\mu=1$, and $c=2$, equation~\eqref{eq:type-one-single} gives
\[
 \Lambda(x)=1-\frac12e^{-x/2}.
\]
For $p_1=p$, $p_2=1-p$, the transform is
\[
 \cL_\Lambda(\theta)=
 \frac{(c-\lambda\mu(2-p))(1+\mu\theta)^2}
 {(c\theta-\lambda)(1+\mu\theta)^2
   +\lambda p(1+\mu\theta)+\lambda(1-p)}.
\]
Its two nonzero poles are the roots
\[
 r_{1,2}=\frac{\lambda\mu-2c\pm
 \sqrt{(\lambda\mu)^2+4\lambda\mu c(1-p)}}{2c\mu},
\]
which are negative under the safety loading condition. The survival function has the form
\[
 \Lambda(x)=1+B e^{r_1x}+C e^{r_2x},
\]
where the coefficients are
\[
 B=\frac{\lambda\mu(2-p)r_1+\lambda(3-2p)}{c(r_2-r_1)},
 \qquad
 C=\frac{\lambda\mu(2-p)r_2+\lambda(3-2p)}{c(r_1-r_2)}.
\]
For $\lambda=\mu=1$, $p=1/2$, and $c=5/2$, this gives
\[
 \Lambda(x)\approx1-0.6266e^{-0.3101x}+0.0266e^{-1.2899x}.
\]

For $p_1=p_2=p_3=1/3$, $\lambda=\mu=1$, and $c=3$, partial-fraction inversion gives
\[
 \Lambda(x)\approx1-0.7017e^{-0.2085x}
 +0.0351e^{-1.2291x}\cos(0.2965x)
 +0.0266e^{-1.2291x}\sin(0.2965x).
\]
The corresponding initial survival probability is $\Lambda(0+)=1/3$.
\end{example}

\begin{figure}[htbp]
 \centering
 \includegraphics[width=.62\textwidth]{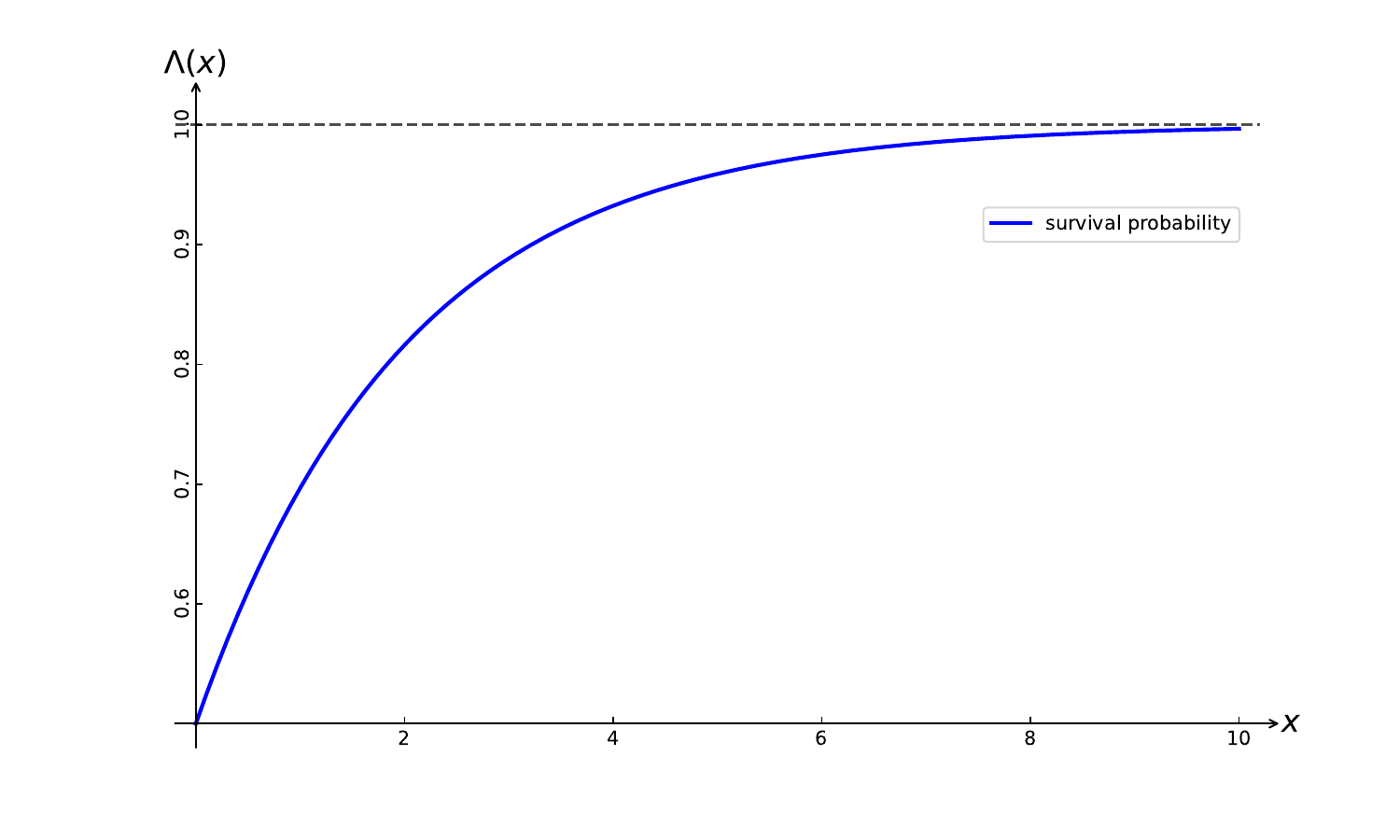}
 \caption{Survival probability in the one-claim-per-event Type I model.}
 \label{fig:type-one-single}
\end{figure}
\begin{figure}[htbp]
 \centering
 \includegraphics[width=.62\textwidth]{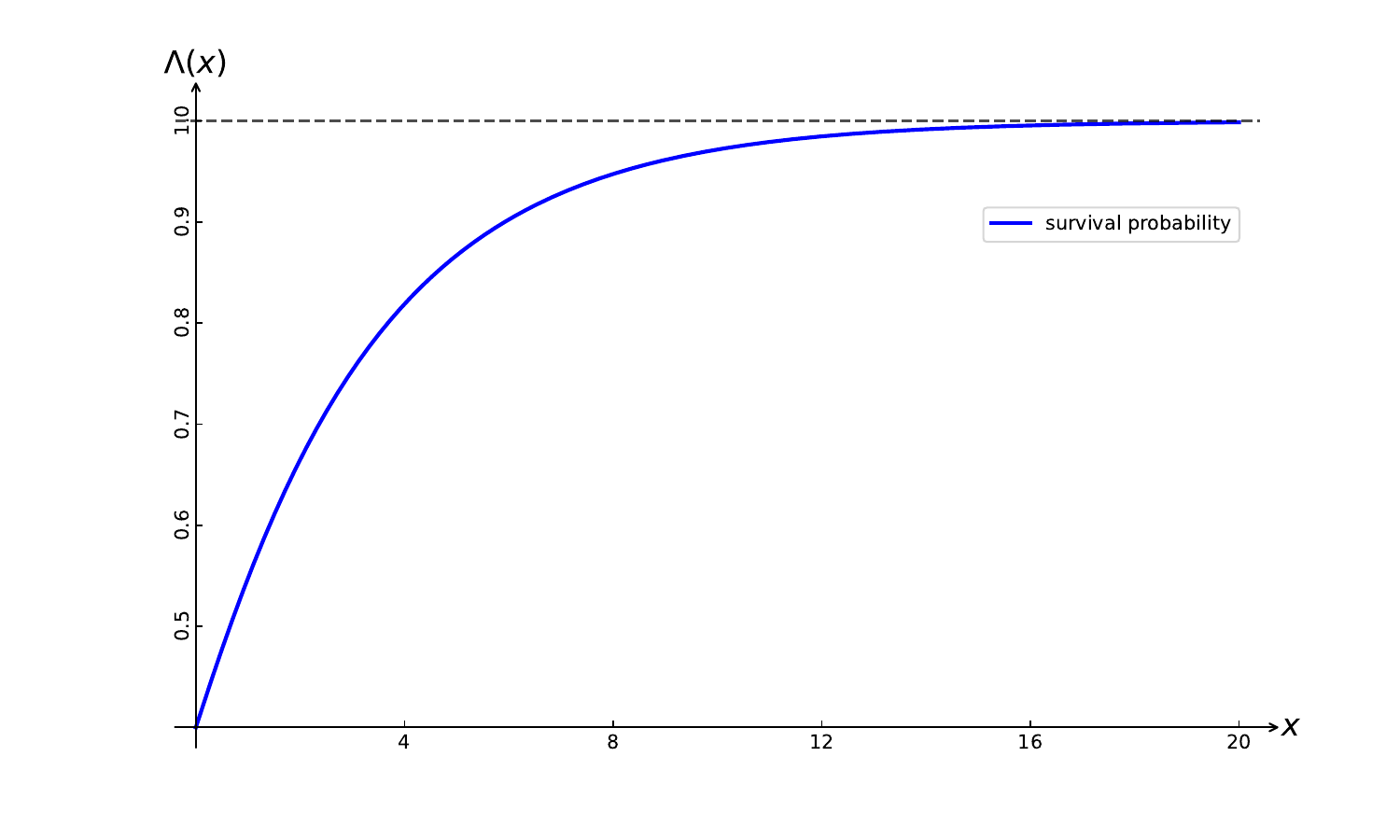}
 \caption{Survival probability when each event produces one or two claims with equal probability.}
 \label{fig:type-one-two}
\end{figure}
\begin{figure}[htbp]
 \centering
 \includegraphics[width=.62\textwidth]{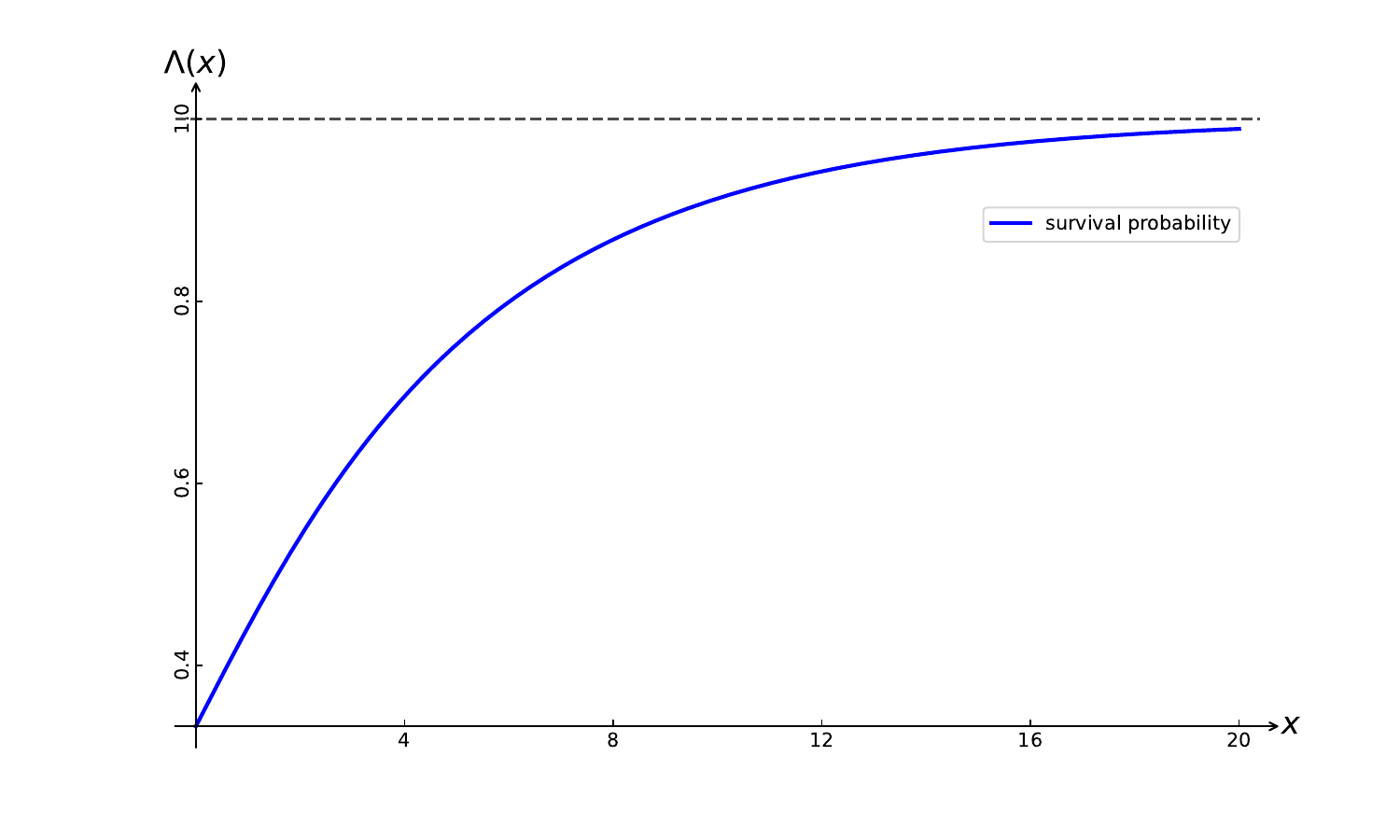}
 \caption{Survival probability when each event produces one, two, or three claims with equal probability.}
 \label{fig:type-one-three}
\end{figure}
\begin{figure}[htbp]
 \centering
 \includegraphics[width=.78\textwidth]{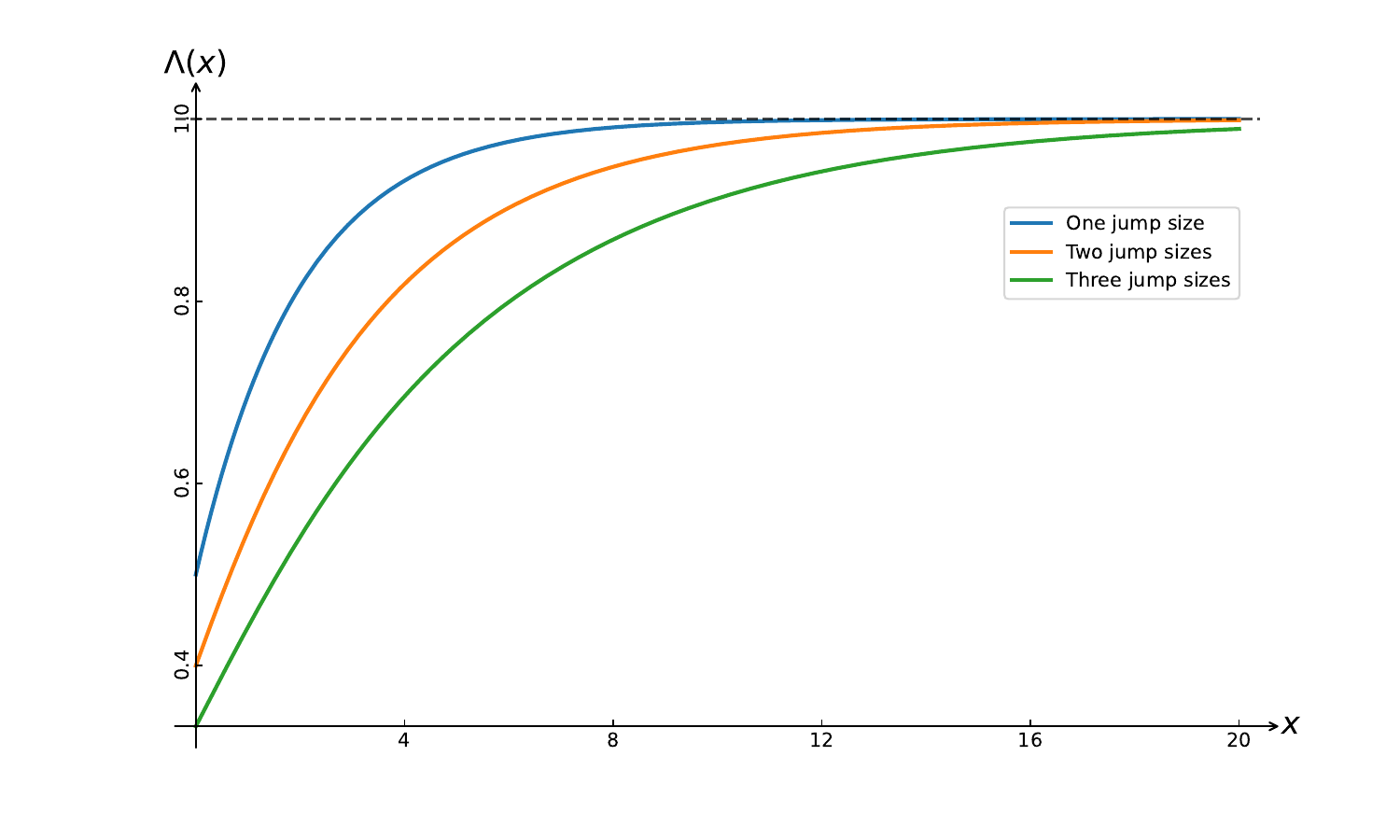}
 \caption{Comparison of the Type I survival probabilities for different claim-cluster distributions.}
 \label{fig:type-one-comparison}
\end{figure}

\begin{remark}
The numerical comparison in \cref{fig:type-one-comparison} is made under parameter choices that keep the mean surplus drift fixed. In these examples, the more dispersed cluster-size specifications have lower survival curves. This is an observation for the selected parameterisations, not a universal stochastic-ordering assertion for all claim-size and cluster-size distributions.
\end{remark}

\subsection{Type II: stochastic premium income and clustered claims}
\label{subsec:type-two}

We now retain the random claim-cluster size but replace the deterministic premium rate by compound-Poisson premium income. Let premium payments arrive according to a Poisson process with rate $\lambdatilde$, with i.i.d. individual premium amounts $\xitilde_i$ having mean $\mutilde$. At each claim-event arrival of a Poisson process with rate $\lambda$, a random number $V$ of individual claims is generated, with
\[
 \Pbb(V=i)=p_i,\qquad i=1,\ldots,n,\qquad \sum_{i=1}^n p_i=1.
\]
The surplus process is
\begin{equation}\label{eq:type-two-surplus}
 X_t(x)=x+\sum_{i=1}^{\Ntilde_t}\xitilde_i
 -\sum_{j=1}^{N_t}\sum_{k=1}^{V_j}\xi_{j,k}.
\end{equation}
The net profit condition is
\begin{equation}\label{eq:type-two-safety}
 \lambdatilde\mutilde>\lambda\mu\sum_{i=1}^n i p_i.
\end{equation}

\begin{theorem}[Laplace transform for Type II]
\label{thm:type-two-transform}
Suppose the individual premium amounts are exponential with mean $\mutilde$. Then
\begin{equation}\label{eq:type-two-transform}
 \boxed{\displaystyle
 \cL_\Lambda(\theta)=
 \frac{\lambdatilde\mutilde-\lambda\mu\sum_{i=1}^n i p_i}
 {(1-\mutilde\theta)
 \left[-\lambda+\lambda\sum_{i=1}^n p_i\hatF(\theta)^i\right]
 +\lambdatilde\mutilde\theta},\qquad\theta>0.}
\end{equation}
If individual claims are also exponential with mean $\mu$, then $\hatF(\theta)=(1+\mu\theta)^{-1}$, and the transform is rational. Its inverse can be obtained from Theorem~\ref{thm:rational-inversion}.
\end{theorem}
\begin{proof}
The first-positive-event equation is
\[
 (\lambda+\lambdatilde)\Lambda(x)
 =\lambda\sum_{i=1}^n p_i\int_{[0,x]}\Lambda(x-y)\,dF^{*i}(y)
 +\lambdatilde\int_{[0,\infty)}\Lambda(x+y)\,d\Ftilde(y).
\]
Taking Laplace transforms, the claim term becomes
\[
 \lambda\cL_\Lambda(\theta)\sum_{i=1}^n p_i\hatF(\theta)^i,
\]
while the exponential premium term becomes
\[
 \frac{\lambdatilde}{1-\mutilde\theta}
 \left[\cL_\Lambda(\theta)-\cL_\Lambda(1/\mutilde)\right].
\]
Rearranging and using $\theta\cL_\Lambda(\theta)\to1$ as $\theta\downarrow0$, together with
\[
 \lim_{\theta\downarrow0}
 \frac{-\lambda+\lambda\sum_{i=1}^n p_i\hatF(\theta)^i}{\theta}
 =-\lambda\mu\sum_{i=1}^n i p_i,
\]
identifies $\lambdatilde\cL_\Lambda(1/\mutilde)$ as the numerator in \eqref{eq:type-two-transform}. Substitution proves the result.
\end{proof}

\begin{example}[Exponential individual claims]
\label{ex:type-two-examples}
When $p_1=1$, the Type II model reduces to a classical two-sided compound-Poisson risk model. With exponential premium and claim amounts,
\[
 \Lambda(x)=1-
 \frac{\lambda(\mu+\mutilde)}{(\lambda+\lambdatilde)\mutilde}
 \exp\left\{-\frac{\lambdatilde\mutilde-\lambda\mu}
 {(\lambda+\lambdatilde)\mu\mutilde}x\right\}.
\]
For $\lambda=\mu=\lambdatilde=1$ and $\mutilde=2$, this becomes
\[
 \Lambda(x)=1-\frac34e^{-x/4}.
\]

For $p_1=p$, $p_2=1-p$, the transform is
\[
 \cL_\Lambda(\theta)=
 \frac{(\lambdatilde\mutilde-\lambda\mu(2-p))(1+\mu\theta)^2}
 {(1-\mutilde\theta)
 [-\lambda(1+\mu\theta)^2+\lambda p(1+\mu\theta)+\lambda(1-p)]
 +\lambdatilde\mutilde\theta(1+\mu\theta)^2}.
\]
For $\lambda=\mu=\lambdatilde=1$, $\mutilde=5/2$, and $p=1/2$, partial-fraction inversion yields
\[
 \Lambda(x)\approx1-0.8187e^{-0.1421x}+0.0187e^{-1.4079x}.
\]
For $p_1=p_2=p_3=1/3$, $\lambda=\mu=\lambdatilde=1$, and $\mutilde=3$, the result is
\[
 \Lambda(x)\approx1-0.8567e^{-0.0916x}
 +0.0234e^{-1.2875x}\cos(0.4024x)
 +0.0179e^{-1.2875x}\sin(0.4024x).
\]
The last two examples both satisfy the safety loading condition.
\end{example}

\begin{figure}[htbp]
 \centering
 \includegraphics[width=.62\textwidth]{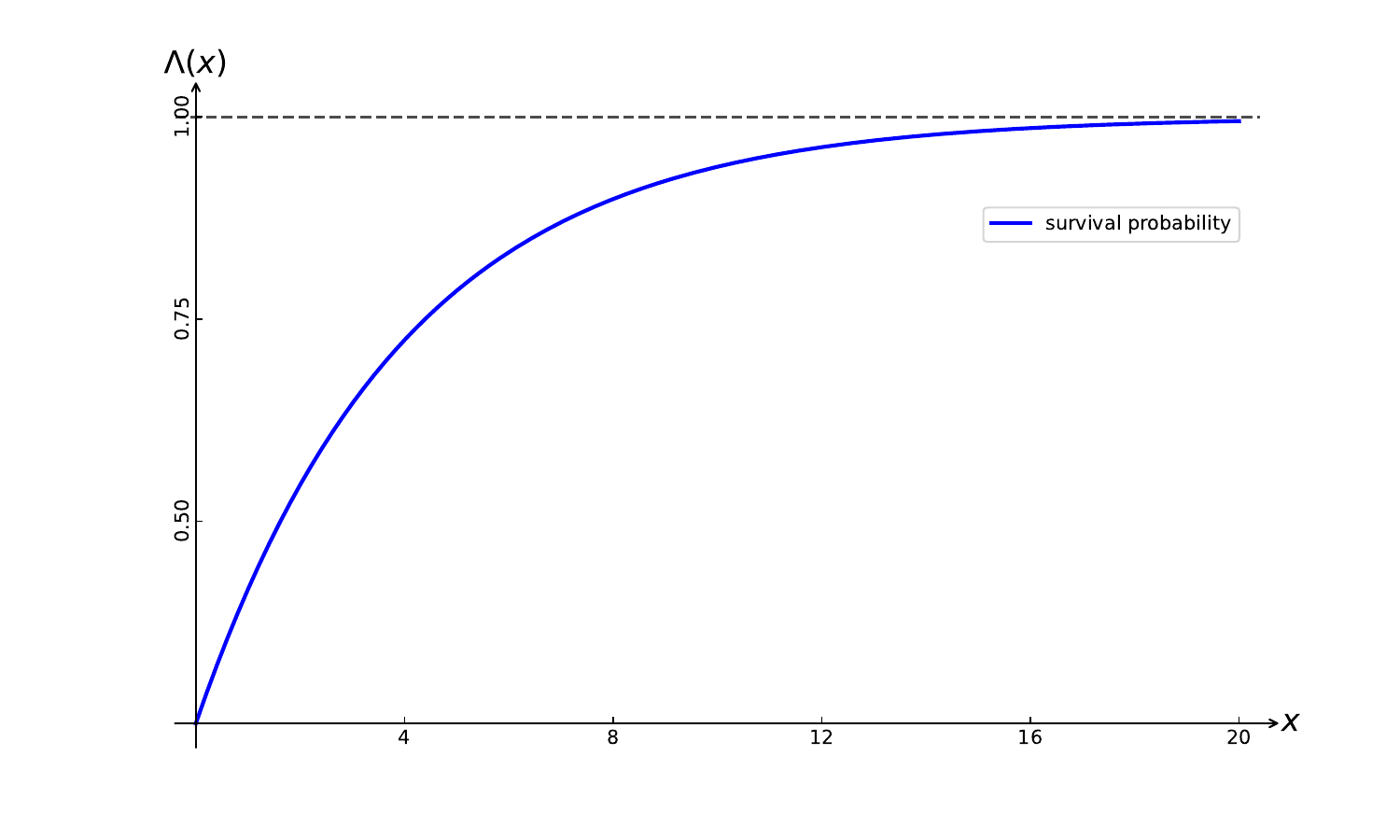}
 \caption{Type II survival probability when each claim event produces one claim.}
 \label{fig:type-two-single}
\end{figure}
\begin{figure}[htbp]
 \centering
 \includegraphics[width=.62\textwidth]{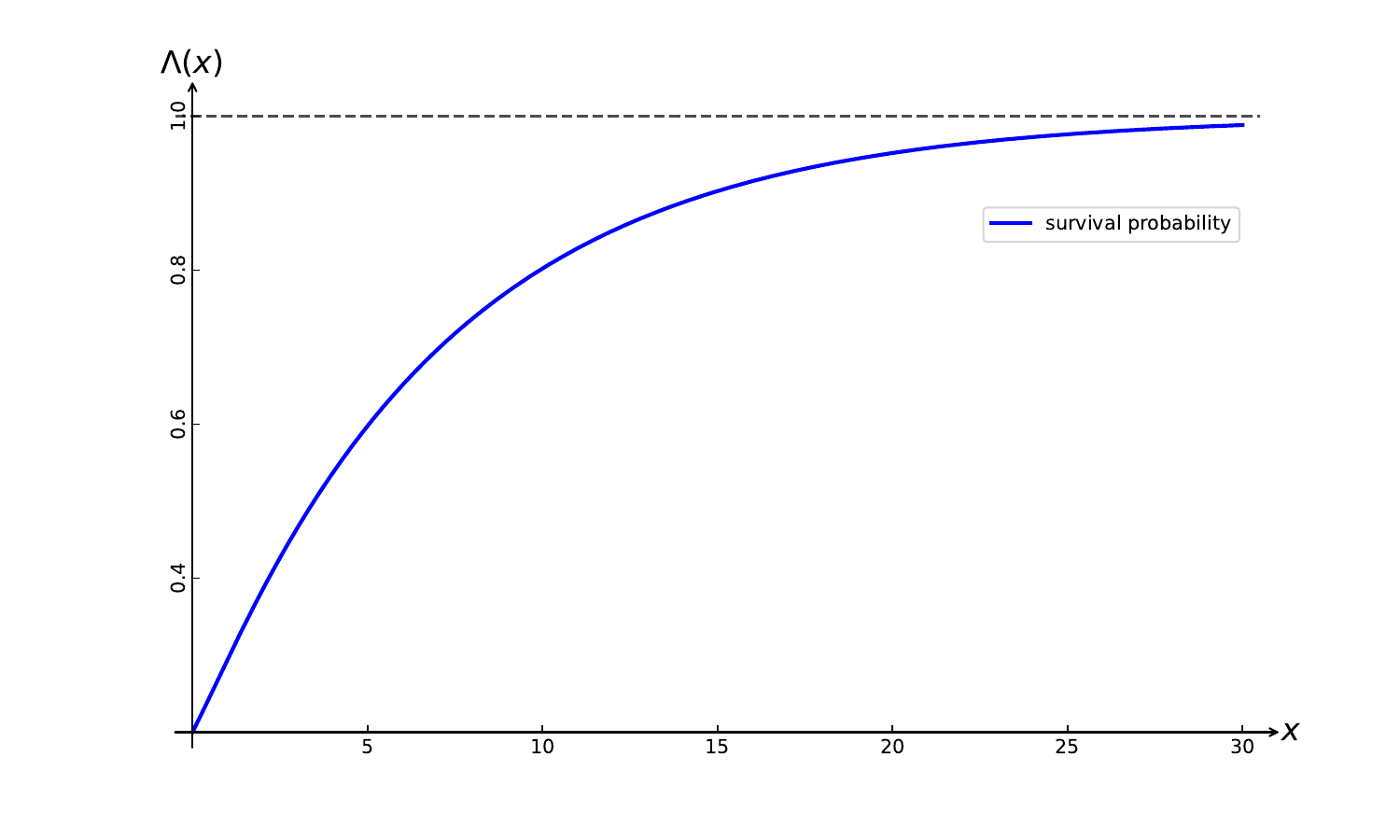}
 \caption{Type II survival probability when each claim event produces one or two claims with equal probability.}
 \label{fig:type-two-two}
\end{figure}
\begin{figure}[htbp]
 \centering
 \includegraphics[width=.62\textwidth]{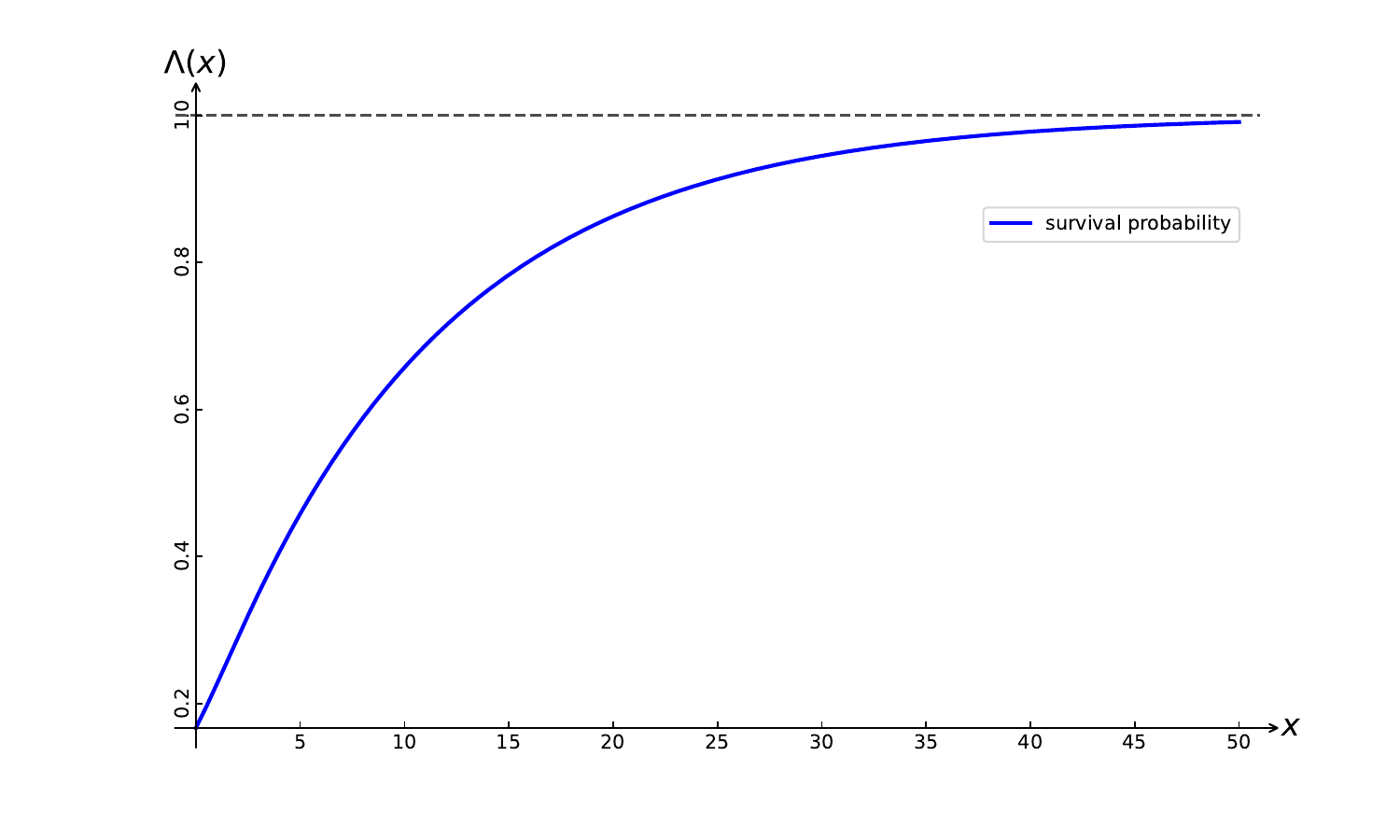}
 \caption{Type II survival probability when each claim event produces one, two, or three claims with equal probability.}
 \label{fig:type-two-three}
\end{figure}
\begin{figure}[htbp]
 \centering
 \includegraphics[width=.78\textwidth]{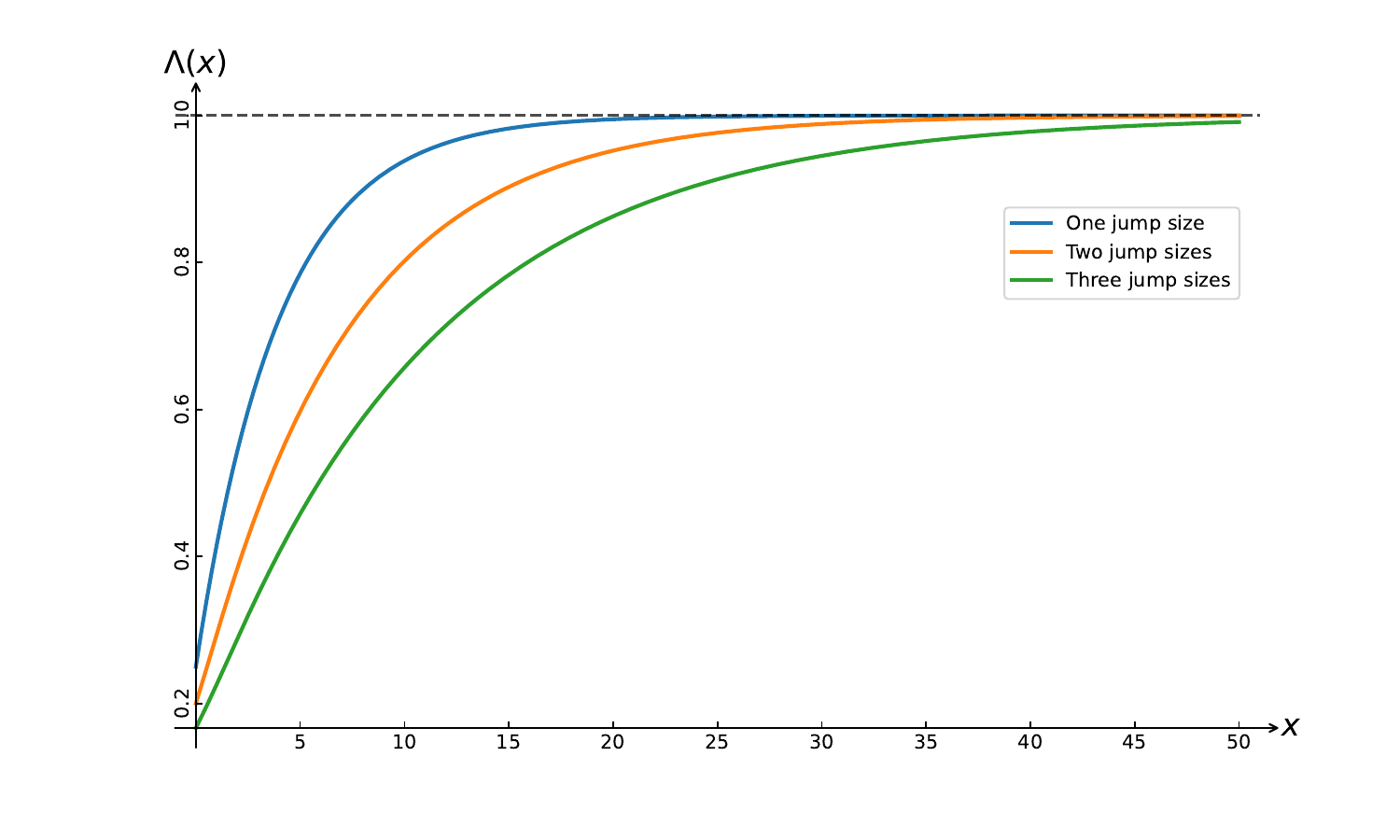}
 \caption{Comparison of Type II survival probabilities under different claim-cluster distributions.}
 \label{fig:type-two-comparison}
\end{figure}

\begin{remark}
For the parameter sets shown in \cref{fig:type-two-comparison}, the curves are compared under matched or closely matched expected premium and claim flows. The numerical results suggest lower survival probabilities for the more dispersed cluster specifications in these examples; no general monotonicity statement is asserted.
\end{remark}

\subsection{Type III: constant and stochastic premium income with clustered claims}
\label{subsec:type-three}

Finally, combine a deterministic premium rate with stochastic premium payments and clustered claims. The surplus process is
\begin{equation}\label{eq:type-three-surplus}
 X_t(x)=x+ct+\sum_{i=1}^{\Ntilde_t}\xitilde_i
 -\sum_{j=1}^{N_t}\sum_{k=1}^{V_j}\xi_{j,k},
\end{equation}
where $N$ and $\Ntilde$ are independent Poisson processes with rates $\lambda$ and $\lambdatilde$, and the cluster size $V$ has probabilities $p_1,\ldots,p_n$. The net profit condition is
\begin{equation}\label{eq:type-three-safety}
 c+\lambdatilde\mutilde>\lambda\mu\sum_{i=1}^n i p_i.
\end{equation}

Suppose that individual claim sizes are exponential with mean $\mu$ and premium amounts are exponential with mean $\mutilde$. Substituting
\[
 \hatF(\theta)=\frac{1}{1+\mu\theta},\qquad
 \hatH(\theta)=\sum_{k=1}^n\frac{p_k}{(1+\mu\theta)^k}
\]
into Theorem~\ref{thm:laplace-constant-stochastic} gives
\begin{equation}\label{eq:type-three-transform}
 \cL_\Lambda(\theta)=
 \frac{c+\lambdatilde\mutilde-\lambda\mu\sum_{i=1}^n i p_i
       -c\mutilde\theta\Lambda(0+)}
 {(c\theta-\lambda+\lambda\hatH(\theta))(1-\mutilde\theta)
       +\lambdatilde\mutilde\theta}.
\end{equation}
After multiplication by $(1+\mu\theta)^n$, the denominator becomes a polynomial. More explicitly, define
\begin{align}
 D_n(\theta)={}&\left[(c\theta-\lambda)(1+\mu\theta)^n
 +\lambda\sum_{k=1}^n p_k(1+\mu\theta)^{n-k}\right](1-\mutilde\theta)\notag\\
 &+\lambdatilde\mutilde\theta(1+\mu\theta)^n.
 \label{eq:type-three-denominator}
\end{align}

\begin{proposition}[Initial survival probability in the rational case]
\label{prop:type-three-initial}
Assume that $D_n$ has a unique positive zero $\theta_1$ and that the pole at $\theta_1$ in \eqref{eq:type-three-transform} is removable. Then
\begin{equation}\label{eq:type-three-initial}
 \boxed{\displaystyle
 \Lambda(0+)=
 \frac{c+\lambdatilde\mutilde-\lambda\mu\sum_{i=1}^n i p_i}
 {c\mutilde\theta_1}.}
\end{equation}
The remaining poles and residues determine $\Lambda(x)$ through Theorem~\ref{thm:rational-inversion}.
\end{proposition}
\begin{proof}
The Laplace transform of a bounded function is finite for every $\theta>0$. Hence the positive zero of the denominator must be cancelled by the numerator. Substituting $\theta=\theta_1$ in the numerator of \eqref{eq:type-three-transform} gives
\[
 c+\lambdatilde\mutilde-\lambda\mu\sum_{i=1}^n i p_i
 -c\mutilde\theta_1\Lambda(0+)=0,
\]
which proves the stated expression.
\end{proof}

\subsubsection{Examples}

\begin{example}[One possible claim per event]
\label{ex:type-three-one}
Take $c=\lambda=\lambdatilde=\mu=\mutilde=1$ and $p_1=1$. Then
\[
 \cL_\Lambda(\theta)=
 \frac{\theta+1}{(\sqrt2+1)\theta(\theta+\sqrt2-1)}
 =\frac1\theta-\frac{2-\sqrt2}{\theta+\sqrt2-1}.
\]
Consequently,
\begin{equation}\label{eq:type-three-one-explicit}
 \Lambda(x)=1-(2-\sqrt2)e^{-(\sqrt2-1)x},\qquad x\geq0,
\end{equation}
and $\Lambda(0+)=\sqrt2-1$.
\end{example}

\begin{figure}[htbp]
 \centering
 \includegraphics[width=.62\textwidth]{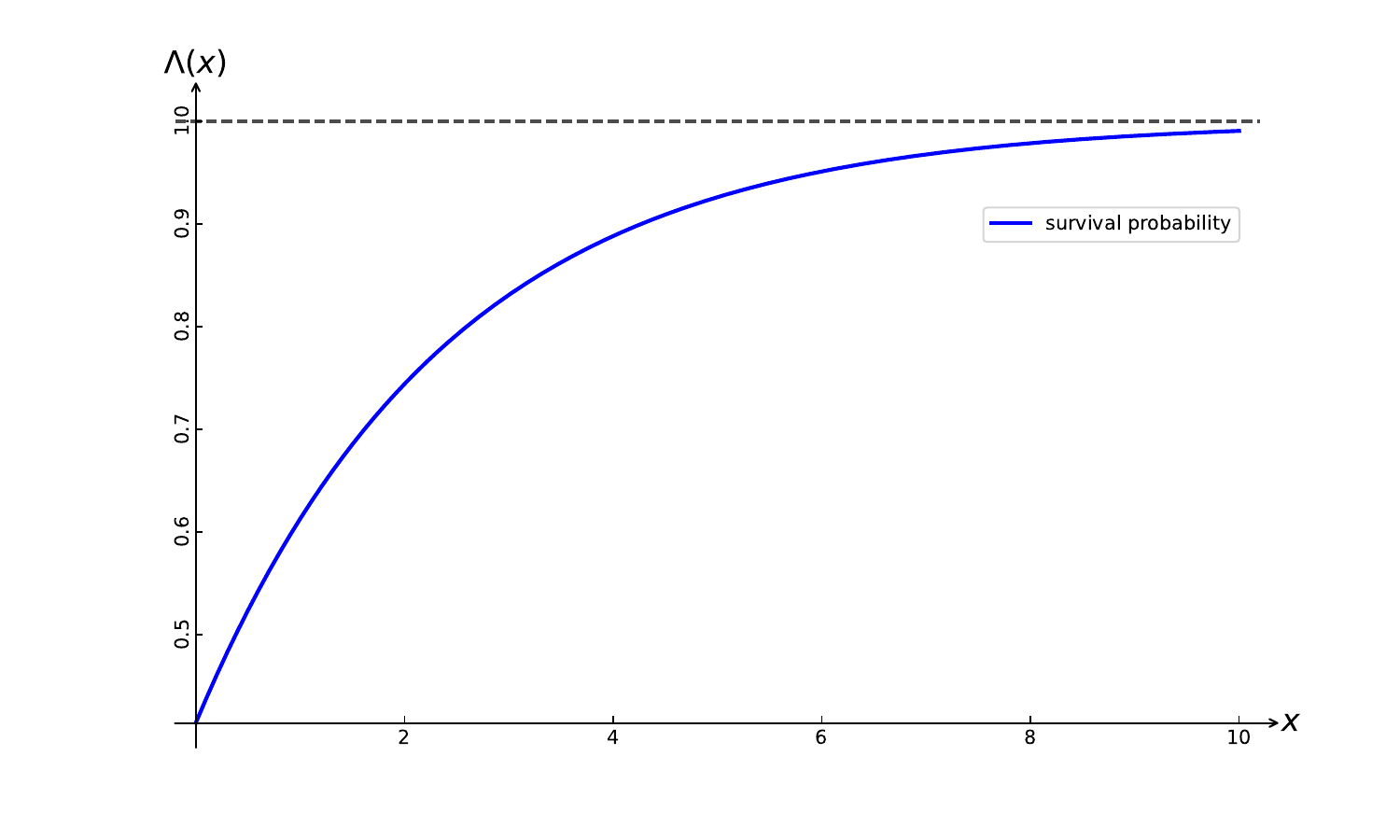}
 \caption{Type III survival probability when each claim event produces one claim.}
 \label{fig:type-three-single}
\end{figure}

\begin{example}[Two possible claims per event]
\label{ex:type-three-two}
Let $c=3/2$, $\lambda=\lambdatilde=\mu=\mutilde=1$, and $p_1=p_2=1/2$. The positive zero of the denominator polynomial is approximately $\theta_1=1.9112$. The other nonzero roots are approximately $-0.2659$ and $-1.3120$, and the initial survival probability is $\Lambda(0+)\approx0.3488$. Partial fractions give
\[
 \cL_\Lambda(\theta)\approx
 \frac1\theta-\frac{0.6759}{\theta+0.2659}
 +\frac{0.0247}{\theta+1.3120},
\]
so that
\begin{equation}\label{eq:type-three-two-explicit}
 \Lambda(x)\approx1-0.6759e^{-0.2659x}+0.0247e^{-1.3120x}.
\end{equation}
\end{example}

\begin{figure}[htbp]
 \centering
 \includegraphics[width=.62\textwidth]{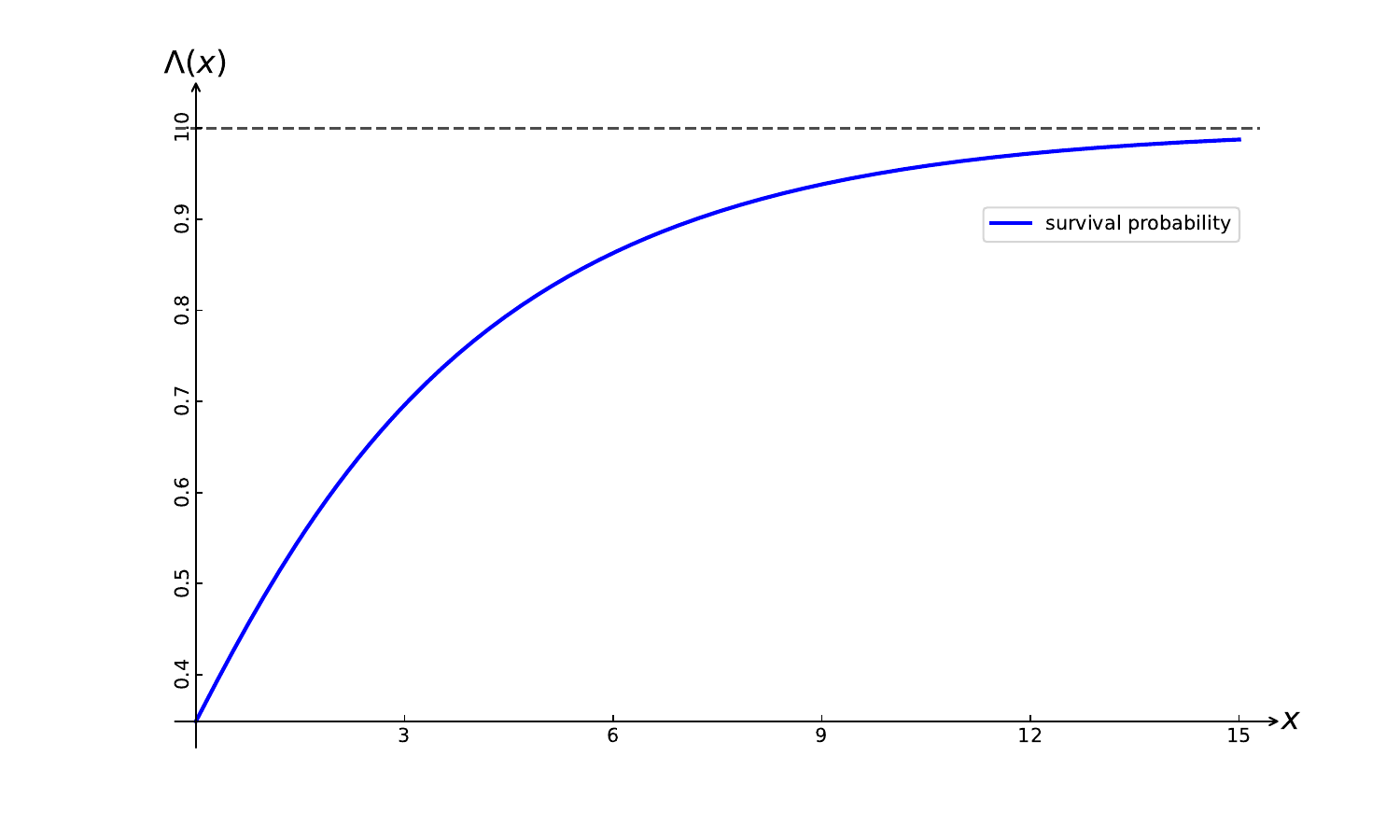}
 \caption{Type III survival probability when the claim-cluster size is equally likely to be one or two.}
 \label{fig:type-three-two}
\end{figure}

\begin{example}[Three possible claims per event]
\label{ex:type-three-three}
Let $c=2$, $\lambda=\lambdatilde=\mu=\mutilde=1$, and $p_1=p_2=p_3=1/3$. The positive denominator zero is approximately $\theta_1=1.6612$; the other poles consist of a negative real pole and a complex-conjugate pair. Inverting the rational transform gives
\begin{align}
 \Lambda(x)\approx{}&1-0.7321e^{-0.1846x}\notag\\
 &+0.0331e^{-1.2383x}\cos(0.3112x)
 +0.0251e^{-1.2383x}\sin(0.3112x),\qquad x\geq0.
 \label{eq:type-three-three-explicit}
\end{align}
The initial survival probability is approximately $0.3010$.
\end{example}

\begin{figure}[htbp]
 \centering
 \includegraphics[width=.62\textwidth]{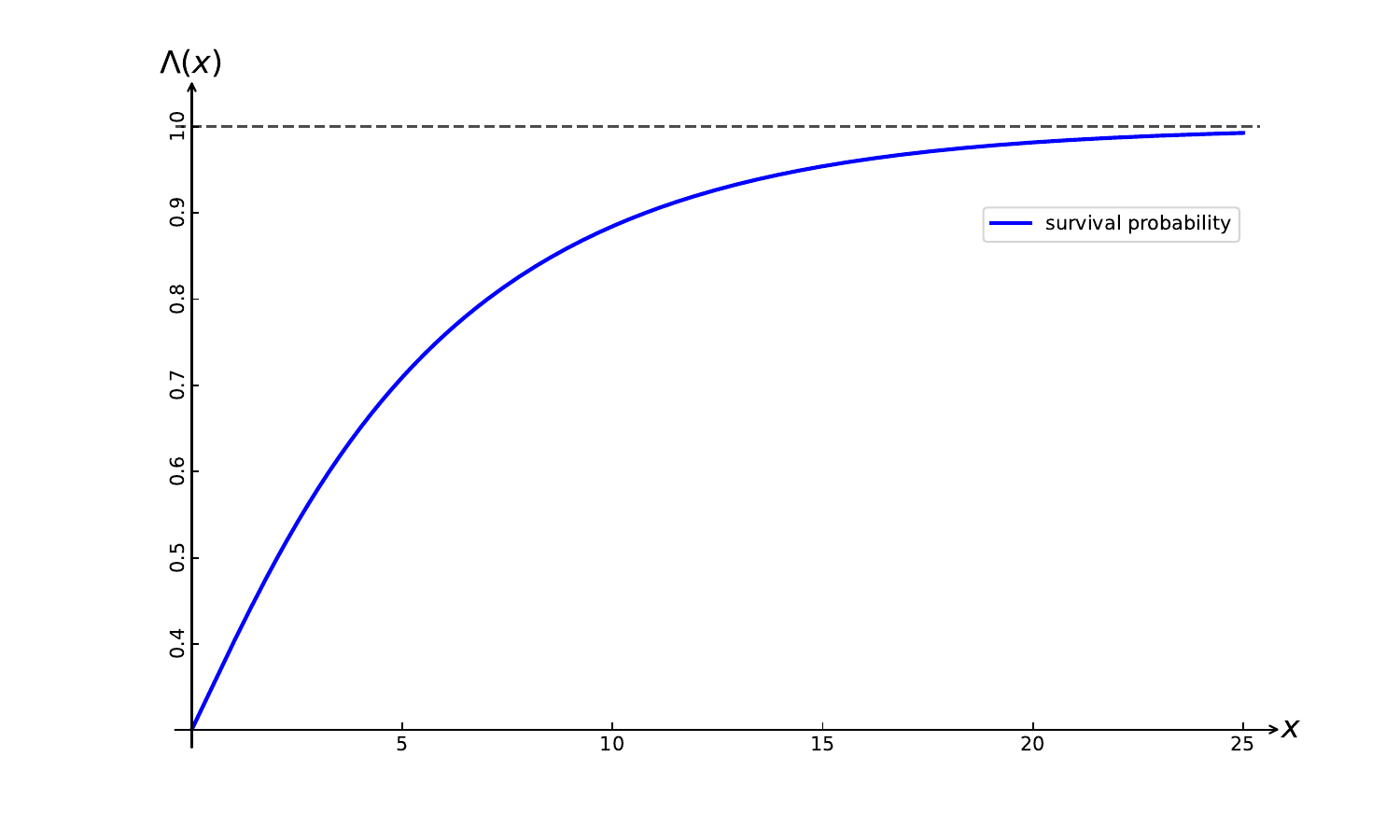}
 \caption{Type III survival probability when the claim-cluster size is equally likely to be one, two, or three.}
 \label{fig:type-three-three}
\end{figure}

\begin{figure}[htbp]
 \centering
 \includegraphics[width=.78\textwidth]{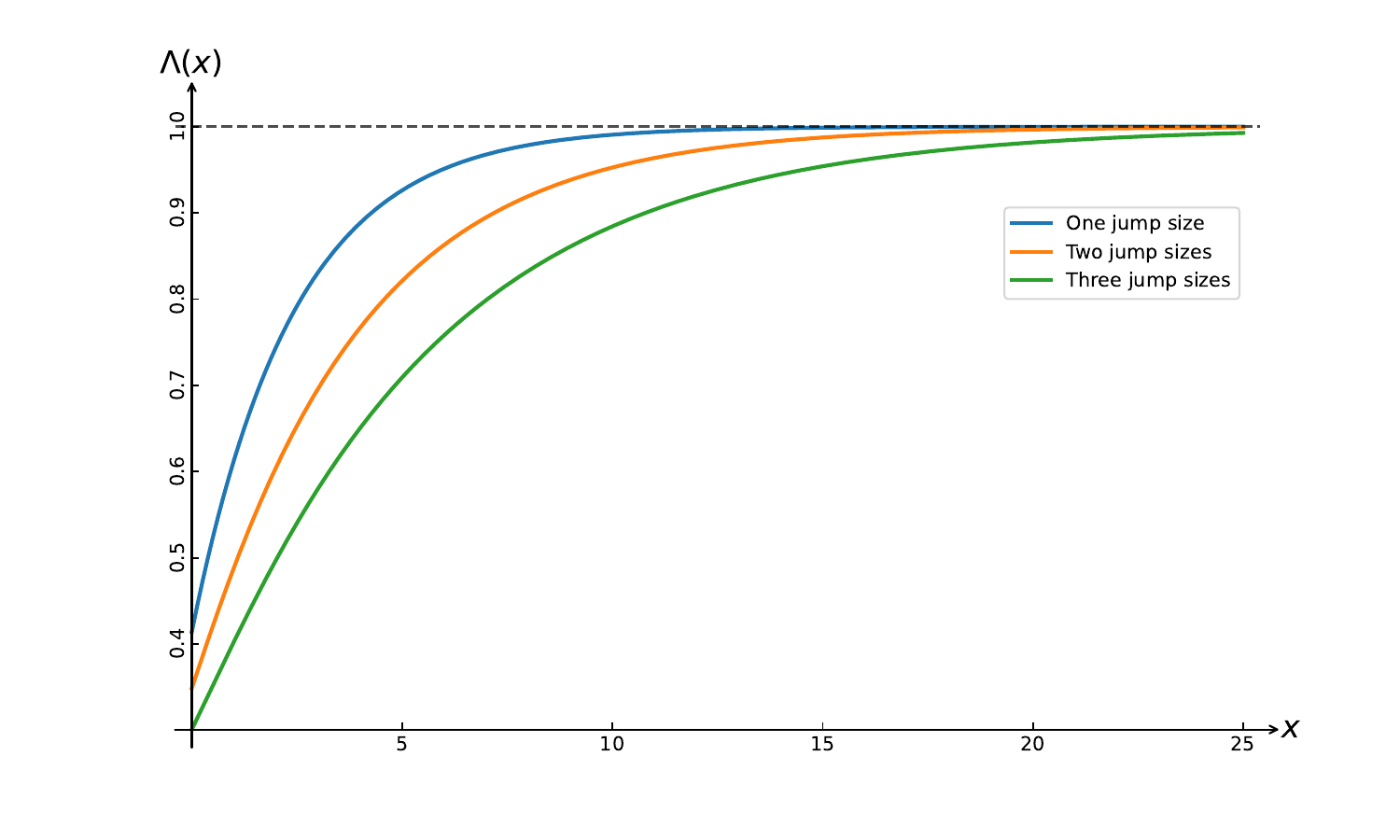}
 \caption{Comparison of Type III survival probabilities for the three claim-cluster distributions.}
 \label{fig:type-three-comparison}
\end{figure}

\begin{remark}
The three numerical illustrations in Type III are calibrated to have the same positive expected net surplus drift. Their survival curves decrease as the chosen cluster-size distribution becomes more dispersed. As in Types I and II, this is a numerical observation for the parameter sets considered rather than a universal ordering theorem.
\end{remark}

\begin{figure}[htb]
\centering
\begin{subfigure}{0.32\textwidth}
    \centering
    \includegraphics[width=\textwidth]{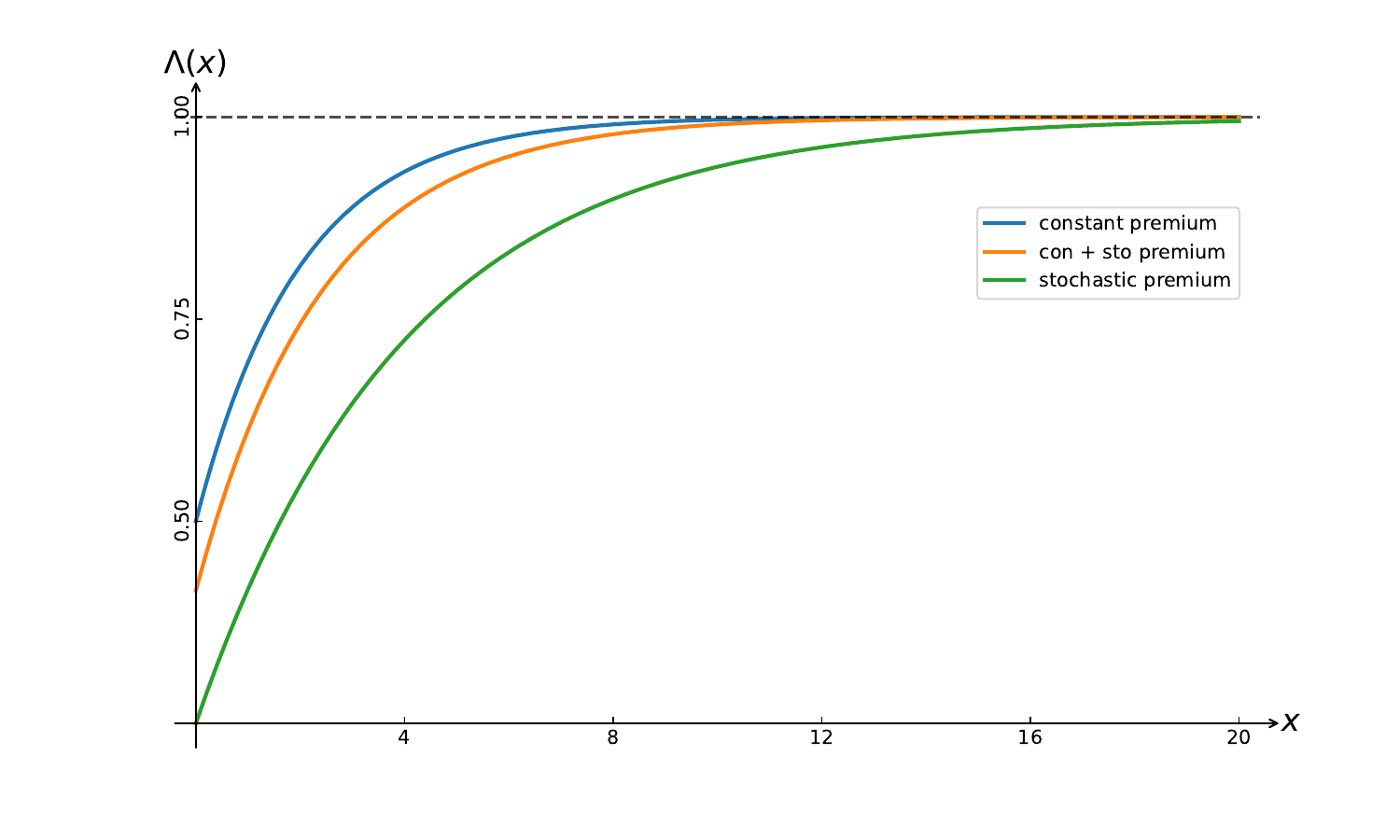}
    \caption{One jump size}
\end{subfigure}
\hfill
\begin{subfigure}{0.32\textwidth}
    \centering
    \includegraphics[width=\textwidth]{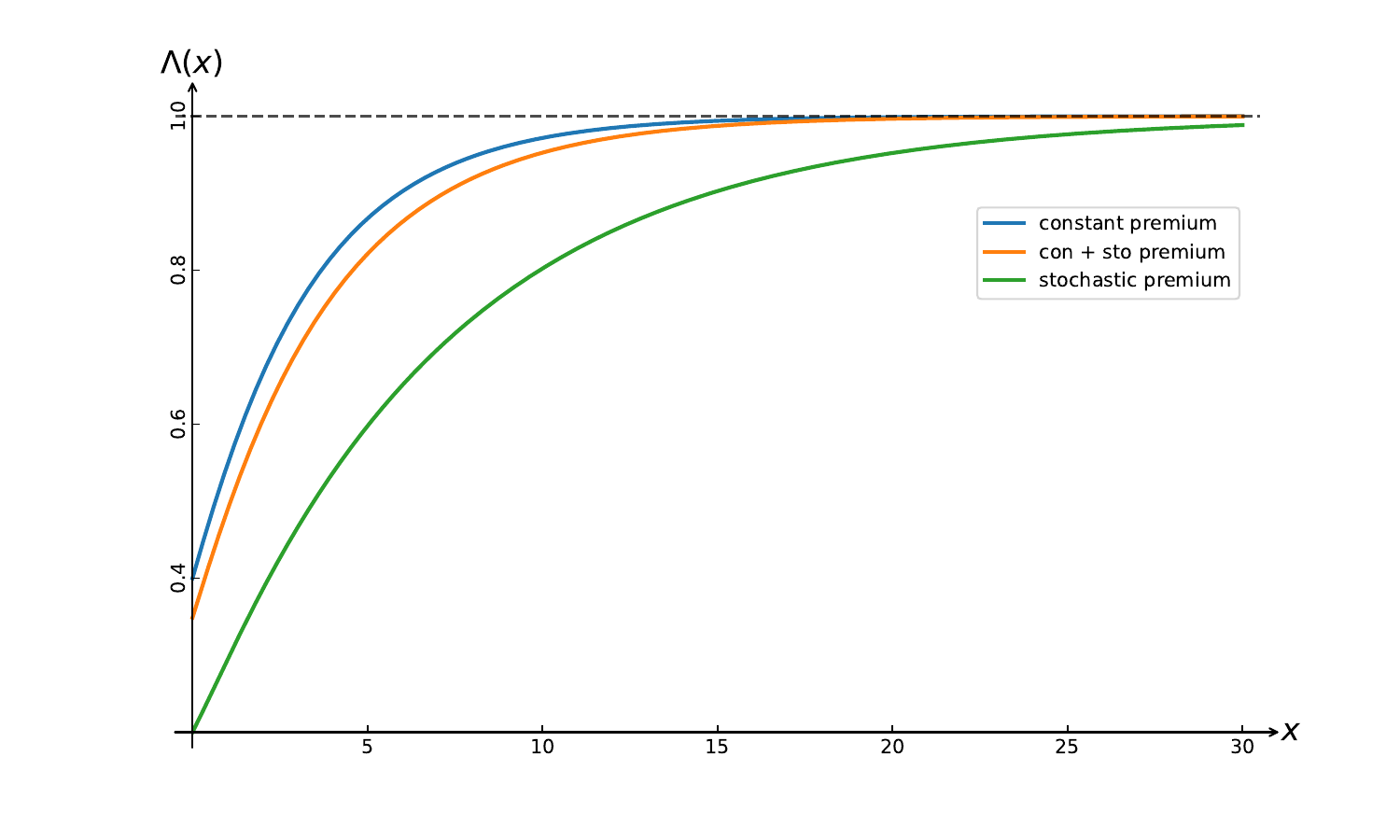}
    \caption{Two jump sizes}
\end{subfigure}
\hfill
\begin{subfigure}{0.32\textwidth}
    \centering
    \includegraphics[width=\textwidth]{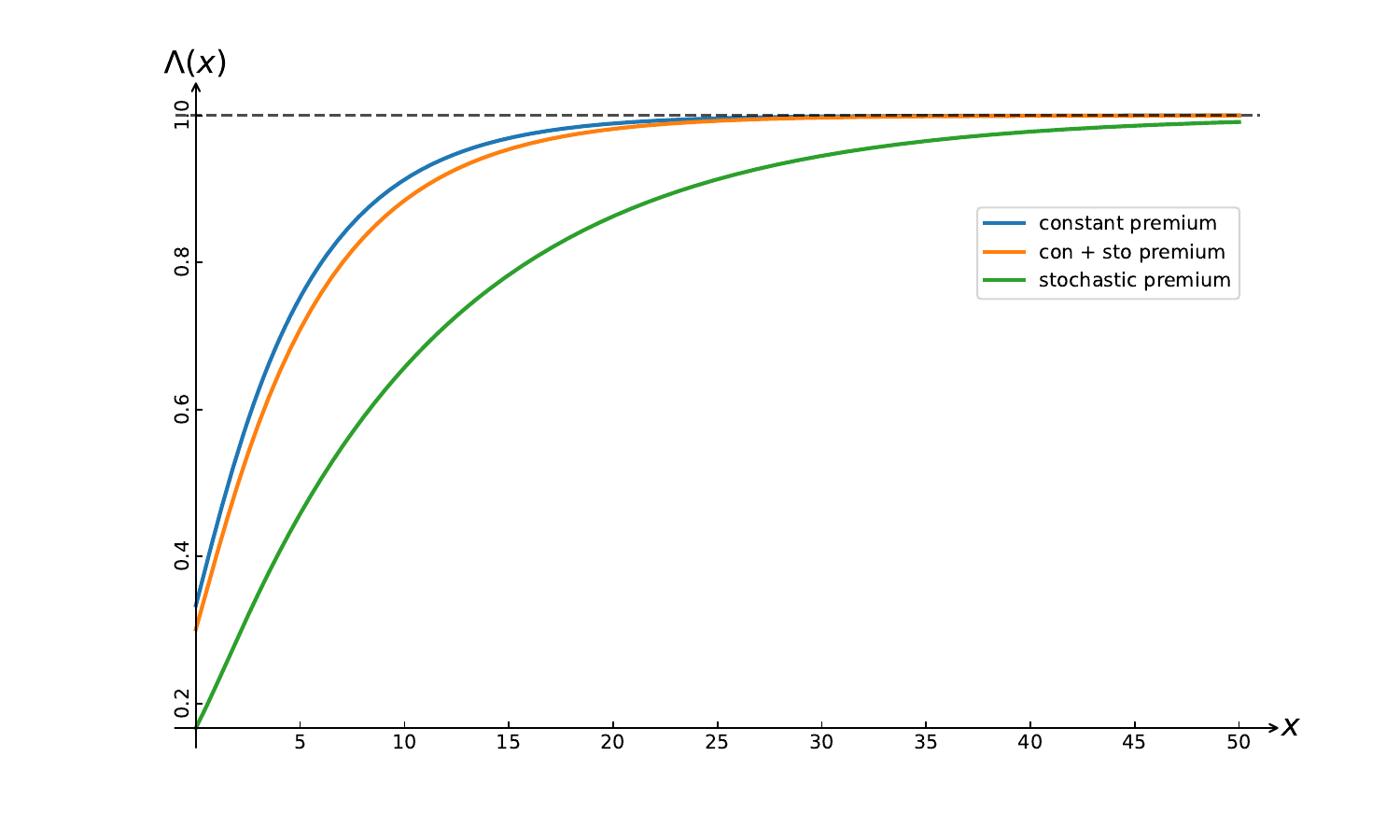}
    \caption{Three jump sizes}
\end{subfigure}
\caption{Comparison of constant, stochastic, and combined premium specifications for one-, two-, and three-claim clusters.}
 \label{fig:all-model-comparison}
\end{figure}

\begin{remark}
Figure~\ref{fig:all-model-comparison} compares the three premium specifications at matched expected net drift within each claim-cluster configuration. In these examples, adding stochasticity to premium income is associated with lower survival probabilities. The figure is intended as a model comparison under the stated calibrations; it does not establish that stochastic premium income always reduces survival for every admissible choice of parameters and distributions.
\end{remark}

\section{Concluding remarks}
\label{sec:conclusion}

This paper has examined ruin models in which premium receipts and claim payments may both be stochastic and clustered. The MIPP formulation provides a common representation for a broad family of compound-Poisson count processes. For the general model, we derived a first-cluster renewal equation and an exponential lower bound for the ultimate survival probability whenever a positive adjustment coefficient exists. The special case of ordinary Poisson premium arrivals admits a tractable Laplace-transform formula when individual premium amounts are exponential. A further compound-Poisson premium specification was developed using an order-two MIPP, with explicit care taken to distinguish individual premium sizes from aggregate cluster sizes.

For models that combine deterministic and stochastic premium income, we obtained an integro-differential equation and a martingale-based ruin bound. When stochastic premium amounts are exponential, the survival transform can be expressed in terms of the claim-size transform and the initial survival probability. If claim amounts are exponential as well, the transform becomes rational and leads to an explicit survival function.

The finite-support examples illustrate how different cluster-size distributions affect survival probabilities under matched mean-drift conditions. The calculations show that cluster dispersion can have a material impact on solvency measures even when expected premium and claim flows are held fixed. A useful next step is to move beyond finite-support cluster sizes and investigate more general dependence structures, including Markov-modulated or state-dependent clustering, as well as empirical calibration of the premium-arrival mechanism.

\end{document}